\documentclass[10pt,a4paper]{article}
\usepackage[margin=2.5cm]{geometry}\usepackage[utf8]{inputenc}
\usepackage[T1]{fontenc}
\usepackage{lmodern}
\usepackage{amsmath,amsfonts,amsthm,amssymb,mathtools}
\usepackage{xcolor,microtype,graphicx,caption}
\usepackage{tikz}
\usetikzlibrary{arrows.meta}
\usepackage{hyperref}
\colorlet{darkblue}{blue!50!black}
\hypersetup{
  colorlinks=true,
  linkcolor=darkblue,
  citecolor=blue,
  urlcolor=blue,
  filecolor=red,
  pdfstartview={FitH}
}
\theoremstyle{plain}
\newtheorem{theorem}{Theorem}[section]
\newtheorem{lemma}[theorem]{Lemma}
\newtheorem{proposition}[theorem]{Proposition}
\newtheorem{corollary}[theorem]{Corollary}
\newtheorem*{maintheorem}{Main Theorem}
\theoremstyle{definition}
\newtheorem{definition}[theorem]{Definition}
\theoremstyle{remark}
\newtheorem{remark}[theorem]{Remark}
\theoremstyle{remark}

\theoremstyle{plain}
 
\newcommand{\R}{\mathbb{R}}
\newcommand{\C}{\mathbb{C}}
\newcommand{\Z}{\mathbb{Z}}

\newcommand{\T}{\mathbb{T}}
\newcommand{\HH}{\mathcal{H}}
\newcommand{\cR}{\mathcal{R}}
\newcommand{\cL}{\mathcal{L}}
\newcommand{\cT}{\mathcal{T}}
\newcommand{\cS}{\mathcal{S}}
\newcommand{\bS}{\mathbb{S}}
\newcommand{\cP}{\mathcal{P}}

\newcommand{\cF}{\mathcal{F}}

\newcommand{\cG}{\mathcal{G}}

\newcommand{\Id}{\mathrm{Id}}
\newcommand{\dd}{\mathrm{d}}
\newcommand{\e}{\varepsilon}
\newcommand{\lan}{\langle}
\newcommand{\ran}{\rangle}
\DeclareMathOperator{\Real}{Re}

\DeclareMathOperator{\dv}{div}
\DeclareMathOperator{\tr}{tr}

\numberwithin{equation}{section}
\numberwithin{figure}{section}
 
\title{Approximate controllability in Sobolev spaces of periodic nonlinear Schr\"odinger equations with bilinear controls}
\author{Vahagn Nersesyan\thanks{NYU-ECNU Institute of Mathematical Sciences at NYU Shanghai, 3663 Zhongshan Road North, Shanghai, 200062, China; e-mail: vahagn.nersesyan@nyu.edu}}
\date{\today}
 
\begin{document}
\maketitle

\begin{abstract}
We consider the nonlinear Schr\"odinger equation
on the torus $\T^d$ with a smooth real-valued potential
and nonlinearity $\kappa|\psi|^{2p}\psi$, where the
dimension $d\ge1$, the integer $p\ge1$, and the
coefficient~$\kappa\in\R$ are arbitrary.
The controls act bilinearly through finitely many low
Fourier modes.~For every integer~$s>d/2$, we prove
approximate controllability in~$H^s$ between any two
nowhere vanishing states on the unit sphere in $L^2$
at any prescribed positive time. The proof combines
small-time phase control, a transport approximation
in~$H^s$, and a construction of short uncontrolled
trajectories connecting states with different
winding~vectors.
\end{abstract}

 

\setcounter{tocdepth}{2}
\tableofcontents

\section*{Introduction}
\addcontentsline{toc}{section}{Introduction}

We consider the nonlinear Schr\"odinger (NLS) equation
on the torus $\T^d=\R^d/2\pi\Z^d$:
\begin{equation}\label{E:0.1}
i\partial_t\psi
=-\Delta\psi+V(x)\psi+\kappa|\psi|^{2p}\psi
+\lan u(t),Q(x)\ran\psi,
\end{equation}
where $d,p\ge1$ are integers, $\kappa\in\R$, and
$\lan\cdot,\cdot\ran$ denotes the Euclidean scalar product.
We write $Q=(Q_1,\dots,Q_q)$ and assume that
$V,Q_1,\dots,Q_q$ are smooth and real-valued.
For any $T>0$, the control is a function
$u\in L^2(0,T;\R^q)$, whose components are the
time-dependent amplitudes of the fixed spatial profiles
$Q_1,\dots,Q_q$.
Equation~\eqref{E:0.1} is supplemented with the initial
condition
\begin{equation}\label{E:0.2}
\psi(0,x)=\psi_0(x),
\end{equation}
where $\psi_0\in H^s=H^s(\T^d;\C)$ and $s>d/2$
is an integer.
We denote the norm in $H^s$ by $\|\cdot\|_s$.
Problem \eqref{E:0.1},~\eqref{E:0.2} is locally well posed
in $H^s$, and the $L^2$ norm of its solutions is conserved.
We therefore consider initial and target states in the
unit sphere
$$
\bS=\{\psi\in L^2(\T^d;\C):\|\psi\|_{L^2}=1\}.
$$
Our aim is to choose the amplitudes $u_j(t)$ so that
the solution is arbitrarily close in $H^s$ to a given
target state at a prescribed positive time.

This paper is a continuation of~\cite{DN}, where
Duca and the author developed a multiplicative version
of the Agrachev--Sarychev method~\cite{AS05,AS06}
to control the periodic phase of solutions
of~\eqref{E:0.1}.
Let us recall the explicit family of potentials
used there.
We denote by $e_1,\dots,e_d$ the standard basis of $\R^d$
and set
\begin{equation}\label{E:0.3}
K=\{e_1,\dots,e_{d-1},(1,\dots,1)\},
\end{equation}
with the convention $K=\{1\}$ when $d=1$.
We assume that
\begin{equation}\label{E:0.4}
\{\mathbf1,\sin\lan k,x\ran,\cos\lan k,x\ran:k\in K\}
\subset\operatorname{span}\{Q_j:j=1,\dots,q\}.
\end{equation}
Under this assumption, it was proved in~\cite{DN}
that, for any $\psi_0\in H^s$ and
$\theta\in C^\infty(\T^d;\R)$, the system can be
steered by a suitable control from $\psi_0$ arbitrarily
close to $e^{i\theta}\psi_0$ in $H^s$ in arbitrarily
small time.
This phase-control result is recalled below in
Theorem~\ref{T:1.6}.

Multiplication by a phase does not change the modulus
of the state.
In the present paper, we show that the same potentials
also allow us to change the modulus and to control
the whole state in $H^s$.
More~precisely, we prove the following result.

\phantomsection\label{T:main}
\begin{maintheorem}
Assume \eqref{E:0.4} and let $s>d/2$ be an integer.
For any nowhere vanishing initial and target states
$\psi_0,\psi_1\in H^s(\T^d;\C)\cap\bS$, any time
$T>0$, and any $\e>0$, there are a control
$u\in L^2(0,T;\R^q)$ and a solution
$\psi\in C([0,T];H^s\cap\bS)$ of \eqref{E:0.1}
such that
$$
\psi(0)=\psi_0,\qquad \|\psi(T)-\psi_1\|_s<\e.
$$
\end{maintheorem}

This theorem holds in every dimension $d\ge1$, for every integer
$p\ge1$, every $\kappa\in\R$, and every smooth real-valued
potential~$V$. In particular, it covers both focusing and
defocusing nonlinearities, without any restriction on their
degree or size. At this level of generality, the uncontrolled
equation is not necessarily globally well posed. The role of
the control is therefore twofold: it keeps the solution defined
on the whole interval $[0,T]$ and steers it close to the target
state at time~$T$.   
The $2d+1$ Fourier potentials in~\eqref{E:0.4} are universal:
they do not depend on $s$, $p$, $\kappa$, $V$, $T$, $\e$, or
the initial and target states. Only the amplitudes~$u_j(t)$
depend on these data. We prove a more general version of this result under
a saturation condition on the control potentials;
see Theorem~\ref{T:3.1}.

 Let us denote by
$$
\phi_l(x)=(2\pi)^{-d/2}e^{i\lan l,x\ran},
\qquad l\in\Z^d,
$$
the normalized exponential eigenfunctions of $-\Delta$
on $\T^d$, which we call \emph{eigenstates}.
In particular, the Main~Theorem allows us to steer any
eigenstate $\phi_l$ arbitrarily close to any other
eigenstate~$\phi_m$ in $H^s$ at any given time.
In the case $V=0$, this was previously known
in~$L^2$; see~\cite[Theorem~B]{DN}.
The passage from $L^2$ to $H^s$ is not straightforward.
Indeed, every nowhere vanishing state $\psi\in H^s$
can be written as
$$
\psi=ae^{i\theta}\phi_l,
\qquad
a,\theta\in H^s(\T^d;\R),\quad a>0,
$$
with a unique $l\in\Z^d$, which we call the
\emph{winding vector} of $\psi$;
see Lemma~\ref{L:A.2}.
Multiplication by $e^{i\eta}$ with~$\eta\in H^s(\T^d;\R)$
does not change this vector. To change the modulus, we show that transport by
suitable diffeomorphisms of $\T^d$ can be approximated
in $H^s$ by controlled solutions.
These diffeomorphisms are isotopic to the identity,
so the transport also preserves the winding vector.

In $L^2$, the winding vector is not an obstruction to
approximation by phase multiplication. Indeed, one can
choose a smooth periodic phase $\eta$ such that
$e^{i\eta}\phi_l$ coincides with $\phi_m$ outside a set
of arbitrarily small measure. In $H^s$, however, the
situation is different. The winding vector is stable
under small perturbations in the uniform norm; for
example, if $l\ne m$, then
$$
\|e^{i\eta}\phi_l-e^{i\theta}\phi_m\|_{L^\infty}
=2(2\pi)^{-d/2},
\qquad \eta,\theta\in H^s(\T^d;\R);
$$
see Remark~\ref{R:1.7}. Since $H^s$ is
continuously embedded in $C^0$, the winding vector is
also stable under small perturbations in $H^s$.
Therefore, phase multiplication and transport alone do
not allow us to approximate in $H^s$ a state by another
one with a different winding vector.
To overcome this obstruction, we need an additional
mechanism, and this is one of the key points of the paper.
We show that the winding vector can be changed by the
uncontrolled dynamics itself: there are short uncontrolled
trajectories whose initial and final states have different
winding vectors; see Lemma~\ref{L:3.2}. Combining
such trajectories with phase control and transport between
states with the same winding vector, we obtain approximate
controllability between any two nowhere vanishing states
in $H^s\cap\bS$.

In the case $d=1$, nowhere vanishing states are dense in
$H^s\cap\bS$, so the Main Theorem remains true for any target
state $\psi_1\in H^s\cap\bS$. For $d\ge2$, this density property fails, since some
zeros persist under small uniform perturbations.
However, in any dimension, the target state can be any
real-valued state in~$H^s\cap\bS$, even if it has zeros;
see Remark~\ref{R:3.3} for more details.

Let us note that the approximate controllability in the Main
Theorem cannot be replaced by exact controllability, since
the control potentials are smooth. Indeed, a solution starting
from a smooth state, for example the constant state $\phi_0$,
remains smooth. On~the other hand, every neighborhood of $\phi_0$
in $H^s\cap\bS$ contains positive states that do not belong
to~$H^{s+1}$, and these states cannot be reached exactly
from~$\phi_0$. This is in the spirit of the well-known
obstruction of Ball, Marsden, and Slemrod~\cite{BMS82} for
bilinear control systems.

In the companion paper~\cite{NLSII}, we show that, on the
circle, any nonzero initial state can be steered to a nowhere
vanishing state in arbitrarily small time. Combined with the
Main Theorem, this gives global approximate controllability
in $H^s\cap\bS$ for every integer $s\ge1$, using only the
potentials $\mathbf1$, $\cos x$, and~$\sin x$. When $V=0$,
by adding two explicit real-valued potentials of finite
regularity, we also obtain global exact controllability in
$H^2(\T;\C)\cap\bS$ when $2p\kappa(2\pi)^{-p}>-4$. Both~results hold at any
given time~$T>0$ and without any nonvanishing assumption.

\subsection*{Previous work}

For linear Schr\"odinger equations, local exact
controllability near the ground state is proved in one
space dimension in~\cite{Bea05} and extended
in~\cite{BL10}. 
Global approximate controllability, including in the
multidimensional case, is obtained by geometric control
methods applied to Galerkin
approximations~\cite{CMSB09,BCCS12,BCS14} and by Lyapunov
techniques~\cite{Ner09,Ner10}. These methods also give
approximation in Sobolev norms: in~\cite{BCS14}, under
suitable regularity assumptions on the operators, and
in~\cite{Ner10}, for generic control potentials. However,
they rely on the linear structure of the equation and, in
general, do not give controllability at a prescribed
positive time.

The question of controllability in small time is of a different
nature. Indeed, a positive minimal time may be necessary even
when the control amplitudes are not restricted: this is the case
for certain models on~$\R^d$ with spatially homogeneous electric
fields~\cite{BCT18}, and further obstructions on Riemannian~manifolds~are~given~in~\cite{BBS21}. Examples of equations that are approximately controllable in
arbitrarily small time are given in~\cite{BPexamples}, where
spatial dilations are used to produce time contractions, so
that small-time controllability follows from a large-time
result.

A different approach is based on the control of the phase in
arbitrarily small time. This was developed in~\cite{DN} for the
polynomial NLS equation on $\T^d$. The key observation is
that a rapidly oscillating phase, followed by a short uncontrolled
evolution, produces a new phase: for any $\psi_0\in H^s$ and
$\varphi\in C^\infty(\T^d;\R)$,
\begin{equation}\label{E:0.5}
e^{-i\varphi/\sqrt\e}\,
\cR_\e\big(e^{i\varphi/\sqrt\e}\psi_0,0\big)
\longrightarrow e^{-i|\nabla\varphi|^2}\psi_0
\quad\text{in }H^s\quad\text{as }\e\to0,
\end{equation}
where $\cR_t(v,0)$ denotes the solution of~\eqref{E:0.1} at time $t$
with initial state $v$ and zero control.
Starting from the phases given by the control potentials and
iterating this construction, one obtains a dense set of phases,
provided that the potentials satisfy a saturation property.
This is the case for the explicit family~\eqref{E:0.4}.
As a result, for every integer $p\ge1$ and every
$\kappa\in\R$, multiplication by $e^{i\theta}$, with
$\theta\in H^s(\T^d;\R)$, can be approximated in $H^s$
in arbitrarily small time.

Small-time control of phases has since been studied in various
other settings. In~\cite{CP23}, it is established in $L^2$ for linear
equations on closed manifolds under a saturation
assumption, with applications to the control of
rotating linear molecules. In~\cite{BLS25}, it is used to
study small-time transitions between eigenfunctions with the same
modulus on manifolds with boundary and on quantum graphs. For~polynomial NLS equations on $\R^d$, small-time controllability
properties with bilinear electromagnetic controls are obtained
in~\cite{DP25}.

Phase control alone does not change the modulus of the state.
For linear equations, it is shown in~\cite{BP} that short
uncontrolled evolutions, inserted between phase multiplications,
approximate transport by diffeomorphisms. Using Moser's
theorem~\cite{Mos65} to transport suitable smooth positive
densities, one can then control the modulus, and the remaining
phase is approximated in $L^2$. Combined with a density argument,
this gives small-time global approximate controllability in $L^2$,
up to a constant phase, for linear equations on $\T^d$ and $\R^d$.
A similar result is obtained in~\cite{Poz26} for linear equations
with magnetic drift.

On a general manifold, constructing a finite family of
control potentials with the saturation property is an
additional difficulty. Such families are constructed in~\cite{Par26} on
compact connected normal homogeneous spaces and in~\cite{BDPNash},
using Nash embeddings, on compact connected Riemannian manifolds
without boundary. Combined with the transport approach
of~\cite{BP}, these constructions give small-time global
approximate controllability in $L^2$ for the linear equation
on such manifolds.

Apart from the phase-control results mentioned above, much less
is known for nonlinear equations. In~this~case, the transport
construction requires estimates along the rapidly oscillating
solutions produced by the phase multiplications. Such estimates
are obtained in~\cite{BCP} for the logarithmic NLS equation, and
they lead to small-time global approximate controllability
in~$L^2$, up to a constant phase. Let us note that, in the
logarithmic case, the equation is globally well posed in $L^2$
in any dimension, and its solutions satisfy a stability estimate
that does not depend on the control. These properties play an important role in the approach
of~\cite{BCP}, but they are not available in the case  of the polynomial equation~considered~here.

For polynomial nonlinearities, the known results are of a
different nature. Local exact controllability near a spatially
constant trajectory is proved in~\cite{BL10} for the
one-dimensional cubic equation with Neumann boundary conditions.
In~\cite{CXZ23}, the defocusing cubic equation on the circle with
two bilinear controls is studied in the semiclassical limit, and
the density and momentum are approximately controlled in small
time between states with positive densities. For additive controls acting on finitely
many Fourier modes, approximate controllability of the defocusing
cubic equation on the two-dimensional torus is proved~in~\cite{Sar12}.

In the present paper, we control the whole complex state
of~\eqref{E:0.1} in $H^s$ between nowhere vanishing states, with a fixed coefficient of the
Laplacian and with polynomial nonlinearities of any degree and
sign. Compared with~\cite{BCP}, this requires a transport
approximation at the regularity of the given initial state and
a mechanism that changes the winding vector. These steps are
described below.

\subsection*{From phase control to state control}

The proof consists of three main steps. First, we show that a
large phase followed by a short uncontrolled evolution
transports the initial density along an optimal transport map.
Combined with a final correction of the phase, this gives
approximate controllability in arbitrarily small time between
states with the same winding vector. Next, we construct short
uncontrolled trajectories that change the winding vector.
Finally, we use phase control once more, now to return
repeatedly near an eigenstate, and reach the target at any
prescribed time. We explain these three steps below.

\medskip\noindent
\emph{Step~1. From phase control to control of the modulus.}
The result of~\cite{DN} allows us to choose freely the spatial
phases used in the construction. For any state
$\psi\in H^s$ and any $\theta\in H^s(\T^d;\R)$, we can
steer $\psi$ arbitrarily close to $e^{i\theta}\psi$ in $H^s$
in arbitrarily small time, using the fixed potentials
in~\eqref{E:0.4}; see Theorem~\ref{T:1.6}.
Indeed, the control potentials directly produce the phases in
their linear span, and short uncontrolled evolutions generate
new phases from them, as in~\eqref{E:0.5}. The saturation property ensures that the
phases obtained in this way are dense. Thus we are not
restricted to the span of the potentials, and we can choose
the phase according to the motion of the density that we want
to produce.

To change the modulus, we use the effect of a phase on the
subsequent evolution of the density. The~difference with the
phase-control mechanism is in the size of the phase. In~\eqref{E:0.5}, a phase of size~$\e^{-1/2}$ is followed by an
uncontrolled evolution of duration~$\e$. On this scale, the
mass has no time to move, and after removing the initial phase,
the limit is multiplication by a new phase. Here we use a phase
of size~$1/\tau$, followed by an uncontrolled evolution of
duration~$\tau$. The phase gradient now gives a velocity of
order~$1/\tau$, which moves the mass a distance of order one
during this short time. At this scale, the transport term
survives in the limit. To obtain the transported state, we
follow the phase through the eikonal equation
$$
\partial_\sigma\varphi+|\nabla\varphi|^2=0,
\qquad \varphi(0)=\frac12g,
$$
and remove its value at the final time.

The idea of producing a transport by phase multiplications and
short uncontrolled evolutions goes back to~\cite{BP}, where it
is used to approximate the action of diffeomorphisms and to
prove small-time approximate controllability of linear
equations. In~\cite{BCP}, it is extended to the logarithmic
Schr\"odinger equation, which required new nonlinear estimates.

The phase is determined by the required change of density.
By the periodic version of Brenier's theorem and the regularity
theory of optimal transport, there is a map $\Id+\nabla g$,
with $g\in H^{s+2}$ and~$I+D^2g>0$, which sends the initial
density to a smooth positive approximation of the target
density; see Lemma~\ref{L:2.5}. An important property of this
map is that
$$
f_\sigma(x)=x+\sigma\nabla g(x),\qquad 0\le\sigma\le1
$$
are diffeomorphisms for all $\sigma\in[0,1]$, since
$$
I+\sigma D^2g=(1-\sigma)I+\sigma(I+D^2g)>0.
$$
These maps are
exactly the limiting characteristics of the solution with
initial phase $g/(2\tau)$. Thus the optimal transport map also
gives the phase that moves the mass, and we only need this
family of maps in the argument.

The main difficulty is to justify this limit in $H^s$ for the
given initial state. In the logarithmic case~\cite{BCP}, there
is a stability estimate independent of the control, which
allows to approximate first the initial state by regular
states. We do not have such an estimate here. Moreover, the
state~$e^{ig/(2\tau)}\psi_0$ may have a very large Sobolev
norm. However, after rescaling time and removing the eikonal
phase, the amplitude $a$ satisfies
$$
\partial_\sigma a+2\cT_{\nabla\varphi}a
=i\tau\big(\Delta a-Va-\kappa|a|^{2p}a\big),
$$
where $\cT_{\nabla\varphi}$ is a first-order transport operator.
Thus the dispersive, potential, and nonlinear terms are
multiplied by~$\tau$.
Then uniform estimates for the amplitude
give existence on the whole interval and convergence in~$H^s$;
see Lemma~\ref{L:2.3} and Proposition~\ref{P:2.1}. Finally,
phase control corrects the remaining periodic phase, and we
obtain approximate controllability in arbitrarily small time
between states with the same winding vector.

\medskip\noindent
\emph{Step~2. Changing the winding vector.}
Phase multiplication and transport do not change the winding
vector. Since $H^s$ is continuously embedded in $C^0$, the
winding vector is also stable under small perturbations
in~$H^s$. Therefore, the passage from $L^2$ to~$H^s$ involves
a topological obstruction, and not only stronger estimates.
To change the winding vector, the solution must pass through
a state with zeros. However, by Step~1, states with the same
winding vector can be approximately reached from each other,
so it is enough to find one short uncontrolled trajectory
whose initial and final states have the required winding
vectors. We then steer the initial state close to the initial
state of this trajectory, let the system evolve without
control, and steer the resulting state close to the target.

We first construct such a trajectory for the free equation,
backwards in time. On the real line, we take a function of
modulus one whose phase makes one full turn. Rescaling and
periodization give a short free trajectory on the circle from
winding number zero to one, and hence on $\T^d$ from $l$ to~$m$.
For the nonlinear equation, we start from the same initial
state. Since the final free state has constant modulus, a small
error in the uniform norm does not change its winding vector,
so it is enough to show that the nonlinear solution stays close
to the free one in this norm. We obtain this using the Wiener
norm, which remains bounded under the rescaling, while the
Sobolev norms may grow; see Lemmas~\ref{L:1.4} and~\ref{L:3.2}.

\medskip\noindent
\emph{Step~3. Reaching the target at a prescribed time.} The constructions of Steps~1 and~2 are carried out in
arbitrarily small time. To reach the target at a given
time~$T$, we must also fill the remaining time. Since the
uncontrolled equation is not necessarily globally well posed,
we cannot simply let the system evolve without control. Instead,
we use phase control once more to return repeatedly near an
eigenstate. Indeed, at an eigenstate $\phi_l$, the first-order
change under the uncontrolled equation is a change of phase,
even when $V$ is not constant:
$$
e^{ih(V+|l|^2+\kappa(2\pi)^{-dp})}
\cR_h(\phi_l,0)=\phi_l+O(h^2)
\quad\text{in }H^s.
$$
This change can be cancelled by a phase multiplication, which
we approximate by controls in arbitrarily small time. Repeating
short uncontrolled evolutions and phase multiplications, we
obtain an approximate return near $\phi_l$ at any prescribed
time. This also keeps the solution defined on the whole
interval, using only local solutions. Combining this with Step~1, we obtain approximate
controllability at any prescribed time between states with the
same winding vector. Together with Step~2, this proves the
Main Theorem.

\medskip
The paper is organized as follows. Section~\ref{S:1} contains
preliminaries on the well-posedness of~\eqref{E:0.1} and on
phase control. In Section~\ref{S:2}, we study the transport
mechanism and prove approximate controllability between states
with the same winding vector. In Section~\ref{S:3}, we change
the winding vector and prove the Main Theorem. Some auxiliary
results are proved in the Appendix.

\section{Preliminaries}\label{S:1}

In this section, we recall the local well-posedness result
for equation~\eqref{E:0.1} and establish some regularity
properties of its solutions. We also compare the
uncontrolled solutions with the free ones over short
time intervals. Next, we discuss the saturation property of the control
potentials and the phase-control result of~\cite{DN}.
From this result, we derive an approximate return
property near the eigenstates. 

Throughout the rest of the paper, we assume that
$V,Q_1,\dots,Q_q$ are smooth and real-valued and that~$s>d/2$~is~an~integer.
 
\subsection{Well-posedness and regularity}\label{SS:1.1}

For~$\psi\in L^2(\T^d;\C)$ and $n\in\Z^d$, let
$$
\hat\psi(n)=(2\pi)^{-d}\int_{\T^d}\psi(x)e^{-i\lan n,x\ran}\dd x,\qquad \lan n\ran=(1+|n|^2)^{1/2}.
$$
For $k\ge0$, the Sobolev space $H^k$ is endowed with the norm
$$
\|\psi\|_k^2=(2\pi)^d\sum_{n\in\Z^d}\langle n\rangle^{2k}|\hat\psi(n)|^2.
$$
We also use the \emph{Wiener algebra}  of functions with absolutely convergent Fourier series:
$$
W(\T^d)=\left\{\psi\in L^2(\T^d;\C):\|\psi\|_{W(\T^d)}:=\sum_{n\in\Z^d}|\hat\psi(n)|<\infty\right\}.
$$
When the dimension is clear from the context, we
simply write $\|\cdot\|_W$.
We shall use the following standard properties;
see \cite[Chapter~I, Section~6]{Katznelson}
and \cite[Lemma~5.5]{CDS}.
\begin{lemma}\label{L:1.1}
There is $C_s>0$ such that
\begin{align}
\|\psi\|_{L^\infty}
&\le\|\psi\|_W\le C_s\|\psi\|_s,
\qquad \psi\in H^s,
\label{E:1.1}\\
\|\psi_1\psi_2\|_W
&\le\|\psi_1\|_W\|\psi_2\|_W,
\qquad \psi_1,\psi_2\in W(\T^d).
\label{E:1.2}
\end{align}
The group $e^{it\Delta}$ is an isometry in both
$H^k$ and $W(\T^d)$ for every~$k\ge0$.
Moreover, there is $C_k>0$ such that
\begin{equation}\label{E:1.3}
\|fg\|_k
\le C_k\bigl(\|f\|_k\|g\|_W+\|f\|_W\|g\|_k\bigr),
\qquad f,g\in H^k\cap W(\T^d).
\end{equation}
\end{lemma}
\begin{proof}
The bounds in~\eqref{E:1.1} follow
from the Fourier series representation and the
Cauchy--Schwarz inequality, since $s>d/2$.
The group $e^{it\Delta}$ preserves the $H^k$ and
Wiener norms because its Fourier multiplier
$e^{-it|n|^2}$ has modulus one.
For the product estimates, we use the convolution
formula
$$
\widehat{fg}(n)
=\sum_{m\in\Z^d}\hat f(m)\hat g(n-m).
$$
Taking absolute values and summing over $n$ gives
\eqref{E:1.2}. To prove~\eqref{E:1.3}, we multiply the
convolution formula by $\lan n\ran^k$ and use
$$
\lan n\ran^k
\le C_k\bigl(\lan m\ran^k+\lan n-m\ran^k\bigr)
$$
to obtain
$$
\lan n\ran^k|\widehat{fg}(n)|
\le C_k\sum_{m\in\Z^d}
\lan m\ran^k|\hat f(m)|\,|\hat g(n-m)|
+C_k\sum_{m\in\Z^d}
|\hat f(m)|\,\lan n-m\ran^k|\hat g(n-m)|.
$$
Taking the $\ell^2$ norm and applying Young's
inequality
$\|a*b\|_{\ell^2}\le\|a\|_{\ell^2}\|b\|_{\ell^1}$
to the convolutions on the right-hand side implies
\eqref{E:1.3}.
\end{proof}
 
We shall use the following local well-posedness
result for equation~\eqref{E:0.1};
see~\cite[Proposition~1.1]{DN}.
 
\begin{proposition}\label{P:1.2}
For any $\psi_0^\star\in H^s$ and
$u^\star\in L^2_{\mathrm{loc}}(\R_+;\R^q)$,
there is a maximal time $T_*>0$ and
a unique solution
$\psi^\star\in C([0,T_*);H^s)$
of \eqref{E:0.1}, \eqref{E:0.2} with
$\psi_0=\psi_0^\star$ and $u=u^\star$.
If $T_*<\infty$, then~$\|\psi^\star(t)\|_s\to\infty$ as~$t\to T_*^-$.
Moreover, for any $0<T<T_*$, there are
$\delta,C>0$ such that, for any $\psi_0\in H^s$
and $u\in L^2(0,T;\R^q)$ satisfying
$$
\|\psi_0-\psi_0^\star\|_s
+\|u-u^\star\|_{L^2(0,T)}<\delta,
$$
the problem \eqref{E:0.1}, \eqref{E:0.2}
has a unique solution $\psi\in C([0,T];H^s)$,
and
$$
\sup_{t\in[0,T]}\|\psi(t)-\psi^\star(t)\|_s
\le C\big(
\|\psi_0-\psi_0^\star\|_s
+\|u-u^\star\|_{L^2(0,T)}
\big).
$$
\end{proposition}
 
For $\psi_0\in H^s$ and $T>0$, we denote by
$\Theta(\psi_0,T)$ the set of controls
$u\in L^2(0,T;\R^q)$ for which the problem
\eqref{E:0.1}, \eqref{E:0.2} has a solution
in $C([0,T];H^s)$.
For $u\in\Theta(\psi_0,T)$, we denote the
corresponding solution at time $t\in[0,T]$
by $\cR_t(\psi_0,u)$.
 
As $\kappa$, $V$, $Q$, and $u$ are real-valued,
the $L^2$ norm of the solution is conserved.
Thus, if $\psi_0\in H^s\cap\bS$ and~$u\in\Theta(\psi_0,T)$, then
$\cR_t(\psi_0,u)\in\bS$ for all $t\in[0,T]$.
In the next lemma, we show that higher Sobolev
regularity is preserved and give a continuation
criterion in terms of the Wiener norm.

\begin{lemma}\label{L:1.3}
Let $k\ge s$ be an integer, $\psi_0\in H^k$, and $u\in L^2_{\mathrm{loc}}(\R_+;\R^q)$. Let $T_*>0$ be the maximal existence time of  $\psi(t)=\cR_t(\psi_0,u)$ in $H^s$. Then $\psi\in C([0,T_*);H^k)$. If $T_*<\infty$, then
$$
\sup_{t<T_*}\|\psi(t)\|_W=\infty.
$$
\end{lemma}
\begin{proof}
Let us denote by $T_k$ the maximal existence time in $H^k$. By uniqueness, the $H^k$ solution agrees with the $H^s$ solution on $[0,T_k)$, and $T_k\le T_*$. The Duhamel formula gives
\begin{equation}\label{E:1.4}
\psi(t)=e^{it\Delta}\psi_0-i\int_0^t e^{i(t-r)\Delta}
\big[V\psi(r)+\kappa|\psi(r)|^{2p}\psi(r)+\lan u(r),Q\ran\psi(r)\big]\,\dd r.
\end{equation}
Iterating~\eqref{E:1.3} and using the algebra
property~\eqref{E:1.2} of $W(\T^d)$, we obtain
$$
\bigl\||\psi|^{2p}\psi\bigr\|_k\le C_{k,p}\|\psi\|_W^{2p}\|\psi\|_k.
$$
Combining \eqref{E:1.4}, the isometry of
$e^{it\Delta}$ in $H^k$,  the boundedness of
multiplication by $V$ and $Q_j$ on $H^k$,
and Gronwall's inequality, we get
$$
\|\psi(t)\|_k\le\|\psi_0\|_k
\exp\left(C_k\int_0^t\big(1+|u(r)|+|\kappa|\|\psi(r)\|_W^{2p}\big)\,\dd r\right), \quad t<T_k.
$$
First, let $k=s$ and assume that $T_*<\infty$. If the Wiener norm remains bounded, the right-hand side is bounded on $[0,T_*)$, since $u\in L^1(0,T_*;\R^q)$. This contradicts the blow-up alternative of Proposition~\ref{P:1.2}.
 
Now let $k>s$. If $T_k<T_*$, continuity of the $H^s$ solution on $[0,T_k]$ and \eqref{E:1.1} give a bound for its Wiener norm on this interval. The preceding estimate then bounds the $H^k$ norm up to $T_k$, again contradicting the blow-up alternative. Thus $T_k=T_*$, which proves preservation of the $H^k$ regularity.
\end{proof}
 
In Section~\ref{SS:3.2}, we shall consider
smooth initial states whose Wiener norms are uniformly
bounded, while their Sobolev norms may be large.
The next lemma gives a common existence time for
the uncontrolled solutions and a uniform bound on
their difference from the free solutions with the
same initial states.
Its~proof follows the argument
of~\cite[Proposition~5.8 and Remark~5.9]{CDS}.
 
\begin{lemma}\label{L:1.4}
For every $N>0$, there are $T_N>0$ and $C_N\ge0$,
depending only on $N$, $V$, $\kappa$, and $p$, with the
following property. If $\psi_a\in H^s$ and
$\|\psi_a\|_W\le N$, then the solution
$\psi(t)=\cR_t(\psi_a,0)$ exists on~$[0,T_N]$ and satisfies
\begin{equation}\label{E:1.5}
\|\psi(t)\|_W\le2N,\qquad
\|\psi(t)-e^{it\Delta}\psi_a\|_W\le C_Nt,
\qquad 0\le t\le T_N.
\end{equation}
\end{lemma} 
\begin{proof}
Let $T_*$ be the maximal existence time of $\psi$ in $H^s$.
By the Duhamel formula, the algebra property \eqref{E:1.2}
of $W(\T^d)$, and the isometry of $e^{it\Delta}$, we have
$$
\|\psi(t)-e^{it\Delta}\psi_a\|_W
\le\int_0^t\big(
\|V\|_W\|\psi(r)\|_W
+|\kappa|\|\psi(r)\|_W^{2p+1}\big)\,\dd r,
\qquad t<T_*.
$$
Let us set
$$
C_N=2N\|V\|_W+|\kappa|(2N)^{2p+1}.
$$
Then, as long as $\|\psi(r)\|_W\le2N$ for $0\le r\le t$,
we obtain
\begin{equation}\label{E:1.6}
\|\psi(t)-e^{it\Delta}\psi_a\|_W\le C_Nt,
\end{equation}
and hence, since $\|e^{it\Delta}\psi_a\|_W=\|\psi_a\|_W\le N$,
$$
\|\psi(t)\|_W\le N+C_Nt.
$$
We now choose
$$
T_N=\frac{N}{2(1+C_N)},
$$
so that $N+C_Nt <2N$ for $t\le T_N$.
By continuity, the Wiener norm of $\psi$ cannot reach
$2N$ on the interval~$[0, T_*\wedge T_N)$. If $T_*\le T_N$, this bound
contradicts the continuation criterion in Lemma~\ref{L:1.3}.
Hence~$T_*>T_N$, and~\eqref{E:1.6} holds on $[0,T_N]$.
\end{proof}
 
\subsection{Control of phases}\label{SS:1.2}
 
We use the following standard notion of approximate controllability.
 
\begin{definition}\label{D:1.5}
Let $\psi_0,\psi_1\in H^s$ and $T>0$. Equation~\eqref{E:0.1} is \emph{approximately controllable} in $H^s$ from
$\psi_0$ to~$\psi_1$ at time $T$ if for every $\e>0$, there is $u\in\Theta(\psi_0,T)$ such that
$$
\|\cR_T(\psi_0,u)-\psi_1\|_s<\e.
$$
\end{definition}
To state the phase-control theorem of~\cite{DN},
we first recall the saturation condition used there.
For $\varphi\in C^{s+2}(\T^d;\R)$, we denote
$$
B(\varphi)=\sum_{j=1}^d(\partial_{x_j}\varphi)^2.
$$
Given a finite-dimensional subspace
$\HH\subset C^{s+2}(\T^d;\R)$, let $\cF(\HH)$
be the largest subspace of $C^{s+2}(\T^d;\R)$
whose elements have the form
$$
\theta_0-\sum_{j=1}^nB(\theta_j),
\qquad n\ge1,\quad \theta_0,\ldots,\theta_n\in\HH.
$$
The space $\cF(\HH)$ is finite-dimensional and
contains $\HH$. We define recursively
$$
\HH_0=\HH,\qquad
\HH_j=\cF(\HH_{j-1}),\quad j\ge1,
$$
and set $\HH_\infty=\bigcup_{j\ge0}\HH_j$.
The subspace $\HH$ is called \emph{saturating}
if $\HH_\infty$ is dense in $C^{s+2}(\T^d;\R)$;
see \cite[Definition~2.1]{DN}.
We shall assume that the following condition is satisfied.
\begin{itemize}
\item[\hypertarget{H1}{\textup{(H$_1$)}}]
The subspace $\HH=\operatorname{span}\{Q_j:\ j=1,\dots,q\}$
is saturating, where $Q=(Q_1,\ldots,Q_q)$ is as in~\eqref{E:0.1}.
\end{itemize}
By~\cite[Proposition~2.6]{DN}, assumption~\eqref{E:0.4}
implies~\textup{(\hyperlink{H1}{H$_1$})} for every integer $s>d/2$.
The result below follows from~\cite[Theorem~2.2]{DN}
by approximating the phase by smooth functions.
\begin{theorem}\label{T:1.6}
Assume \textup{(\hyperlink{H1}{H$_1$})}. For any $\psi_0\in H^s$ and $\theta\in H^s(\T^d;\R)$, there are positive times $\tau_n\to0$ and controls $u_n\in\Theta(\psi_0,\tau_n)$ such that
$$
\cR_{\tau_n}(\psi_0,u_n)\longrightarrow e^{i\theta}\psi_0
\quad\text{in }H^s.
$$
\end{theorem}
The following remark shows that phase control
alone does not imply approximate controllability
in $H^s$ between different eigenstates.
\begin{remark}\label{R:1.7}
For any
$\eta,\theta\in H^s(\T^d;\R)$ and  $l,m\in\Z^d$ with $l\ne m$, we have
\begin{equation}\label{E:1.7}
    \|e^{i\eta}\phi_l-e^{i\theta}\phi_m\|_{L^\infty}
=2(2\pi)^{-d/2}.
\end{equation}
Indeed, both states have modulus $(2\pi)^{-d/2}$, so to prove \eqref{E:1.7} it
suffices to find a point where the states are opposite.
Let $j$ be such that $l_j\ne m_j$. At $x=te_j$,
where $e_j$ is the $j$th coordinate vector, the
phase difference~is
$$
A(t)=\eta(te_j)-\theta(te_j)+(l_j-m_j)t,
\qquad 0\le t\le2\pi.
$$
By periodicity,
$|A(2\pi)-A(0)|=2\pi|l_j-m_j|\ge2\pi$.
By the intermediate value theorem, there is
a point at which $A$ equals an odd multiple of $\pi$.
The two states are opposite at this point, which proves~\eqref{E:1.7}.
 
In particular, $\phi_m$ cannot be approximated
uniformly by changing only the periodic phase of
$\phi_l$. The~same holds in $H^s$ by Sobolev embedding.
\end{remark}
From now until the end of Section~\ref{S:3},
we assume that Condition~\textup{(\hyperlink{H1}{H$_1$})} is verified.

As a consequence of Theorem~\ref{T:1.6}, we obtain the
following approximate return property near eigenstates
for any smooth real-valued potential~$V$. We use this
property below to pass from controllability in arbitrarily
small time to controllability at any prescribed positive time.

\begin{corollary}\label{C:1.8}
Assume \textup{(\hyperlink{H1}{H$_1$})}.
Let $z=\phi_l$, where $l\in\Z^d$. For every $T_0>0$ and $\e>0$, there is $\delta>0$ such that, for every $T\in(0,T_0]$ and $v\in H^s$ with $\|v-z\|_s<\delta$, there is $u\in\Theta(v,T)$ satisfying
$$
\|\cR_T(v,u)-z\|_s<\e.
$$
\end{corollary}
The proof of this corollary is given in Appendix~\ref{SS:A.1}.

\section{Transport and amplitude control}\label{S:2}

Let us fix an integer $s>d/2$.
In this section, we prove approximate controllability
in $H^s$ between states~$ae^{i\theta}\phi_l\in H^s\cap\bS$
with $a>0$ and the same winding vector $l$.
We first show that transport along gradient maps
can be approximated by controlled trajectories.
We then combine periodic optimal transport
with phase control to change both the amplitude
and the periodic phase.

\subsection{Transport along gradient maps}\label{SS:2.1}
 
Let us denote by $\cG$ the set of functions
$g\in H^{s+2}(\T^d;\R)$ such that the matrix
$I+D^2g(x)$ is positive definite for every~$x\in\T^d$,
where $D^2g$ is the Hessian of~$g$.
By Sobolev embedding, these functions belong to $C^2$.
Extending $g$ periodically to $\R^d$, we see that the function $x\mapsto |x|^2/2+g(x)$ is strictly convex.

For $g\in\cG$ and $\sigma\in[0,1]$, let us define
$$
\begin{aligned}
F_\sigma:\R^d&\longrightarrow\R^d,
&F_\sigma(x)&=x+\sigma\nabla g(x),\\
f_\sigma:\T^d&\longrightarrow\T^d,
&f_\sigma([x])&=[F_\sigma(x)],
\end{aligned}
$$
where $[x]$ denotes the class of $x$ modulo $2\pi\Z^d$.
The identity
\begin{equation}\label{E:2.1}
F_\sigma(x+2\pi n)=F_\sigma(x)+2\pi n,
\qquad n\in\Z^d
\end{equation}
shows that $f_\sigma$ is well defined. We write $x$ instead of $[x]$ when no confusion is
possible. By Lemma~\ref{L:A.1},
$F_\sigma$~and $f_\sigma$ are $C^1$ diffeomorphisms
of $\R^d$ and $\T^d$, respectively. The Jacobian
determinant of $f_\sigma$ is
$$
m_\sigma(x)=\det Df_\sigma(x)
=\det\big(I+\sigma D^2g(x)\big)>0.
$$
Following \cite[Section~2.1]{BP}, we associate with
$f_1$ the \emph{transport operator}
$\cL_g:L^2\to L^2$ defined by
\begin{equation}\label{E:2.2}
(\cL_g\psi)(x)
=m_1(f_1^{-1}(x))^{-1/2}\psi(f_1^{-1}(x)).
\end{equation}
The change of variables $y=f_1(x)$ shows that
$\cL_g$ is unitary on $L^2$. By Lemma~\ref{L:A.1}, $\cL_g$ is also
a bounded linear operator on $H^s$.
The following proposition shows that $\cL_g\psi_0$
can be approximated in $H^s$ by controlled solutions
in arbitrarily small time.
\begin{proposition}\label{P:2.1}
Assume \textup{(\hyperlink{H1}{H$_1$})}.
For any $g\in\cG$ and $\psi_0\in H^s$, there are
positive times $t_n\to0$ and controls
$u_n\in\Theta(\psi_0,t_n)$ such that
$$
\cR_{t_n}(\psi_0,u_n)\longrightarrow\cL_g\psi_0
\quad\text{in }H^s.
$$
\end{proposition}
For the linear equation~\eqref{E:0.1}, that is, when
$\kappa=0$, the $L^2$ version of this proposition follows
from the proof of~\cite[Theorem~3]{BP}. For the logarithmic
Schr\"odinger equation, approximation of gradient flows
in $L^2$ is proved in~\cite[Section~4]{BCP}.
Here we approximate the action of $\cL_g$ in $H^s$
for polynomial nonlinearities of arbitrary degree.

We adapt the phase conjugation and uniform Sobolev estimates
from the proof of~\cite[Proposition~1.2]{DN}. At the scale
used there, the limit is a phase multiplication. At the
scale used here, the transport term survives in the limit,
and we use the full solution of the eikonal equation,
as in~\cite[Proposition~2.5]{BCP}.

The proof of Proposition~\ref{P:2.1} is given in
Section~\ref{SS:2.2}. As a preparation, we write the
solutions of the eikonal and transport equations in terms
of the characteristics $x\mapsto f_\sigma(x)$. For $\xi\in H^{s+1}(\T^d;\R^d)$, let~us~set
$$
\cT_\xi\psi=\lan\xi,\nabla\psi\ran
+\frac12(\dv\xi)\psi.
$$
Integration by parts gives
$\Real\lan\cT_\xi\psi,\psi\ran_{L^2}=0$
for every $\psi\in H^1$.

\begin{lemma}\label{L:2.2}
Let $g\in\cG$ and $\psi_0\in H^s$, and define
$\varphi$ and $b$ by
$$
\begin{aligned}
\varphi(\sigma,x)
&=\frac12g(f_\sigma^{-1}(x))
+\frac{\sigma}{4}|\nabla g(f_\sigma^{-1}(x))|^2,\\
b(\sigma,x)
&=m_\sigma(f_\sigma^{-1}(x))^{-1/2}
\psi_0(f_\sigma^{-1}(x)).
\end{aligned}
$$
Then $\varphi\in C([0,1];H^{s+2})\cap C^1([0,1];H^{s+1})$
and $b\in C([0,1];H^s)\cap C^1([0,1];H^{s-1})$.
Moreover, $b(1)=\cL_g\psi_0$, and
\begin{alignat}{2}
\partial_\sigma\varphi+|\nabla\varphi|^2
&=0,
&\qquad \varphi(0)&=\frac12g,
\label{E:2.3}\\
\partial_\sigma b+2\cT_{\nabla\varphi(\sigma)}b
&=0,
&\qquad b(0)&=\psi_0.
\label{E:2.4}
\end{alignat}
\end{lemma}
 
\begin{proof}
Let us write $P_\sigma=F_\sigma^{-1}$.
By Lemma~\ref{L:A.1}, $P_\sigma-\Id$
and $b(\sigma)=\cL_{\sigma g}\psi_0$ depend continuously
on $\sigma$ in~$H^{s+1}$ and $H^s$, respectively.
The formula for $\varphi$ also gives
$\varphi\in C([0,1];H^{s+1})$.
Since $g\in C^2$, this formula can be differentiated
in $x$ and $\sigma$.
Evaluating it at $f_\sigma(x)$, differentiating in $x$,
and using the symmetry of~$Df_\sigma$, we~obtain
$$
Df_\sigma(x)\nabla\varphi(\sigma,f_\sigma(x))
=\frac12Df_\sigma(x)\nabla g(x).
$$
Since $Df_\sigma$ is invertible,
\begin{equation}\label{E:2.5}
\nabla\varphi(\sigma,f_\sigma(x))
=\frac12\nabla g(x).
\end{equation}
Thus $\nabla\varphi=\frac12\nabla g\circ P_\sigma$
belongs to $C([0,1];H^{s+1})$, and hence
$\varphi\in C([0,1];H^{s+2})$.
Differentiating the formula for
$\varphi(\sigma,f_\sigma(x))$ in $\sigma$
and using~\eqref{E:2.5}, we obtain
$$
\partial_\sigma\varphi(\sigma,f_\sigma(x))
=-\frac14|\nabla g(x)|^2
=-|\nabla\varphi(\sigma,f_\sigma(x))|^2.
$$
This proves~\eqref{E:2.3} and shows that
$\varphi\in C^1([0,1];H^{s+1})$.
 
Next we verify the equation for $b$.
Differentiating~\eqref{E:2.5} in $x$
   and using Jacobi's formula, we obtain
$$
\frac{\partial_\sigma m_\sigma(x)}{m_\sigma(x)}
=\tr\big(D^2g(x)(I+\sigma D^2g(x))^{-1}\big)
=2\Delta\varphi(\sigma,f_\sigma(x)).
$$
For any $\chi\in C^\infty(\T^d)$, a change of variables
gives
$$
\int_{\T^d}b(\sigma,y)\chi(y)\,\dd y
=\int_{\T^d}m_\sigma(x)^{1/2}\psi_0(x)
\chi(f_\sigma(x))\,\dd x.
$$
Differentiating in $\sigma$ and using
$\partial_\sigma f_\sigma=2\nabla\varphi(\sigma)\circ f_\sigma$,
we obtain
$$
\frac{\dd}{\dd\sigma}\int_{\T^d}b(\sigma,y)\chi(y)\,\dd y
=\int_{\T^d}b(\sigma,y)
\bigl(2\lan\nabla\varphi(\sigma,y),\nabla\chi(y)\ran
+\Delta\varphi(\sigma,y)\chi(y)\bigr)\,\dd y.
$$
This proves~\eqref{E:2.4}.
Since $\cT_{\nabla\varphi}b\in C([0,1];H^{s-1})$,
it also gives $b\in C^1([0,1];H^{s-1})$.
The initial conditions and the identity
$b(1)=\cL_g\psi_0$ follow from the definitions.
\end{proof}

\subsection{Small-time transport approximation}\label{SS:2.2}

For $g\in\cG$, let us consider the solution
of~\eqref{E:0.1} with $u=0$ and initial state
$e^{ig/(2\tau)}\psi_0$.
Following the proof of~\cite[Proposition~2.5]{BCP},
we rescale time and remove the phase factor
$e^{i\varphi/\tau}$.
To compare the resulting function with $b$ in $H^s$,
we adapt the energy estimates from the proof
of~\cite[Proposition~1.2]{DN} and obtain bounds
uniform in $\tau$.

\begin{lemma}\label{L:2.3}
Let $g\in\cG$ and $\psi_0\in H^s$, and let
$\varphi$ and $b$ be as in Lemma~\ref{L:2.2}.
For all sufficiently small $\tau>0$, the solution
$\cR_t(e^{ig/(2\tau)}\psi_0,0)$ is defined on $[0,\tau]$ and
$$
\sup_{0\le\sigma\le1}
\big\|e^{-i\varphi(\sigma)/\tau}
\cR_{\tau\sigma}(e^{ig/(2\tau)}\psi_0,0)-b(\sigma)\big\|_s
\longrightarrow0
\qquad\text{as }\tau\to0.
$$
\end{lemma}

\begin{proof}
\emph{Step~1. Uniform in $\tau$ estimates.}
For $h\in H^{s+2}$ and $\tau\in(0,1]$, let us set
$$
a(\sigma)=a_\tau[h](\sigma)
=e^{-i\varphi(\sigma)/\tau}
\cR_{\tau\sigma}(e^{ig/(2\tau)}h,0)
$$
on its interval of existence in $[0,1]$.
By Lemmas~\ref{L:1.3} and~\ref{L:2.2}, the map
$\sigma\mapsto a_\tau[h](\sigma)$ is continuous
with values in $H^{s+2}$ and continuously
differentiable with values in $H^s$.
Substitution into~\eqref{E:0.1} and the eikonal
equation~give
\begin{equation}\label{E:2.6}
\partial_\sigma a+2\cT_{\nabla\varphi(\sigma)}a
=i\tau\big(\Delta a-Va-\kappa|a|^{2p}a\big),
\qquad a(0)=h.
\end{equation}
We use the norm
$$
|a|_s^2=\sum_{|\alpha|\le s}
\|\partial^\alpha a\|_{L^2}^2,
$$
which is equivalent to $\|\cdot\|_s$.
By the Leibniz rule and Sobolev embedding,
$$
\big\|[\partial^\alpha,\cT_{\nabla\varphi(\sigma)}]a
\big\|_{L^2}
\le C_s\|\varphi(\sigma)\|_{s+2}|a|_s,
\qquad |\alpha|\le s.
$$
Since $\varphi$ is bounded in $H^{s+2}$,
integration by parts and the algebra property
of $H^s$ imply
$$
\frac12\frac{\dd}{\dd\sigma}|a|_s^2
\le C|a|_s^2+C\tau|a|_s^{2p+2},
$$
where $C\ge1$ is independent of $\tau$ and $h$.
Let us fix $R\ge1$, assume that $\|h\|_s\le R$, and denote
$$
M_R=2e^{C+1}R,\qquad
\tau_R=\frac1{CM_R^{2p}}.
$$
For $0<\tau\le\tau_R$, as long as $|a|_s\le M_R$,
the preceding inequality gives
$$
\frac12\frac{\dd}{\dd\sigma}|a|_s^2
\le(C+1)|a|_s^2.
$$
Since $|h|_s\le\|h\|_s\le R$, Gronwall's inequality
implies that
$$
|a(\sigma)|_s
\le e^{(C+1)\sigma}|h|_s
\le e^{C+1}R=M_R/2,\qquad \sigma\le1.
$$
By continuity, $|a|_s$ cannot reach $M_R$ on its
interval of existence in $[0,1]$.
 
For each $0<\tau\le\tau_R$, multiplication by
$e^{i\varphi(\sigma)/\tau}$ is bounded on $H^s$,
uniformly in $\sigma\in[0,1]$. Therefore,
the continuation in
Proposition~\ref{P:1.2}   shows that
$\cR_t(e^{ig/(2\tau)}h,0)$ is defined on $[0,\tau]$.
By~equivalence of norms, we   have
\begin{equation}\label{E:2.7}
\sup_{0\le\sigma\le1}\|a_\tau[h](\sigma)\|_s\le C_R,
\qquad 0<\tau\le\tau_R,
\end{equation}
where $C_R$ is independent of $\tau$ and $h$
with $\|h\|_s\le R$.
 
Applying the same energy estimate to the difference
of two solutions and using the local Lipschitz
continuity of $a\mapsto|a|^{2p}a$ in $H^s$, we obtain,
after increasing $C_R$ if necessary,
\begin{equation}\label{E:2.8}
\sup_{0\le\sigma\le1}
\|a_\tau[h](\sigma)-a_\tau[\tilde h](\sigma)\|_s
\le C_R\|h-\tilde h\|_s,
\qquad 0<\tau\le\tau_R
\end{equation}
for $h,\tilde h\in H^{s+2}$ with
$\|h\|_s,\|\tilde h\|_s\le R$.
 
Let us now take $h\in H^s$ with $\|h\|_s\le R$
and choose a sequence $h_n\in H^{s+2}$ such that
$h_n\to h$ in $H^s$ and $\|h_n\|_s\le R$.
For each fixed $0<\tau\le\tau_R$,
estimate~\eqref{E:2.8} implies that
$a_\tau[h_n]$ converges in $C([0,1];H^s)$.
Multiplying by $e^{i\varphi(\sigma)/\tau}$ and
returning to the time variable $t=\tau\sigma$,
we see that the corresponding solutions
of~\eqref{E:0.1} converge in $C([0,\tau];H^s)$.
Passing to the limit in the Duhamel
formula~\eqref{E:1.4}, we find that the limit
is a solution with initial state $e^{ig/(2\tau)}h$
and $u=0$.
By uniqueness, it coincides with
$\cR_t(e^{ig/(2\tau)}h,0)$.
Thus $a_\tau[h_n]\to a_\tau[h]$ in $C([0,1];H^s)$,
and the existence result and the preceding estimates
remain valid for initial states in $H^s$.
 
\emph{Step~2. Convergence to $b$ in $H^s$.}
Let us fix $R\ge\max\{1,\|\psi_0\|_s\}$ and take
$0<\tau\le\tau_R$, with $\tau_R$ as in Step~1.
Since $b$ need not belong to $H^{s+2}$, we first
regularize it in space. Let $J_\delta$, $0<\delta\le1$,
be a standard convolution regularization on $\T^d$,
and set
$$
b_\delta=J_\delta b,\qquad
r_\delta=2[\cT_{\nabla\varphi},J_\delta]b.
$$
By the transport equation~\eqref{E:2.4},
$$
\partial_\sigma b_\delta
+2\cT_{\nabla\varphi}b_\delta=r_\delta.
$$
For $\xi\in H^{s+1}(\T^d;\R^d)$ and $h\in H^s$,
integration by parts, the Leibniz rule, and
Sobolev embedding give
\begin{equation}\label{E:2.9}
\|[\cT_\xi,J_\delta]h\|_s
\le C_s\|\xi\|_{s+1}\|h\|_s,
\end{equation}
with $C_s$ independent of $\delta\in(0,1]$.
For smooth $\xi,h$, the commutator converges to zero
in $H^s$ as $\delta\to0$.
The bound~\eqref{E:2.9} and density extend this
convergence to all $\xi\in H^{s+1}$ and $h\in H^s$,
uniformly on compact subsets of $H^{s+1}\times H^s$.
Since $\nabla\varphi\in C([0,1];H^{s+1})$
and $b\in C([0,1];H^s)$, it follows that
$$
\sup_{0\le\sigma\le1}\|b_\delta(\sigma)-b(\sigma)\|_s
+\sup_{0\le\sigma\le1}\|r_\delta(\sigma)\|_s
\longrightarrow0
\qquad\text{as }\delta\to0.
$$
 
Let us set $a=a_\tau[\psi_0]$ and
$v_\delta=a-b_\delta$.
Subtracting the equation for $b_\delta$
from~\eqref{E:2.6}, we obtain
\begin{equation}\label{E:2.10}
\begin{aligned}
\partial_\sigma v_\delta
+2\cT_{\nabla\varphi(\sigma)}v_\delta
&=i\tau(\Delta v_\delta-Vv_\delta) +i\tau(\Delta b_\delta-Vb_\delta-\kappa|a|^{2p}a)
-r_\delta,\\
v_\delta(0)&=\psi_0-J_\delta\psi_0.
\end{aligned}
\end{equation}
For fixed $\delta$, we have
$b_\delta\in C([0,1];H^{s+2})$.
Let us choose $h_n\in H^{s+2}$ converging to $\psi_0$ in
$H^s$, with~$\|h_n\|_s\le R$, and set
$v_{n,\delta}=a_\tau[h_n]-b_\delta$.
Repeating the energy calculation of Step~1 and
using~\eqref{E:2.7}, we obtain
$$
\frac12\frac{\dd}{\dd\sigma}|v_{n,\delta}|_s^2
\le C|v_{n,\delta}|_s^2
+C\bigl(\tau C_\delta+\|r_\delta\|_s\bigr)
|v_{n,\delta}|_s.
$$
Here $C$ is independent of $n,\tau,\delta$, while
$C_\delta$ may depend on $\delta$ but not on
$n$ or $\tau$.
Since $v_{n,\delta}(0)=h_n-J_\delta\psi_0$,
Gronwall's inequality gives
$$
\sup_{0\le\sigma\le1}\|v_{n,\delta}(\sigma)\|_s
\le C\left(
\|h_n-J_\delta\psi_0\|_s
+\int_0^1\|r_\delta(\sigma)\|_s\,\dd\sigma
+\tau C_\delta\right).
$$
Passing to the limit by~\eqref{E:2.8}
and using $b(0)=\psi_0$, we obtain
$$
\sup_{0\le\sigma\le1}\|a(\sigma)-b(\sigma)\|_s
\le C\left(
\sup_{0\le\sigma\le1}\|b_\delta(\sigma)-b(\sigma)\|_s
+\int_0^1\|r_\delta(\sigma)\|_s\,\dd\sigma
+\tau C_\delta\right).
$$
Letting first $\tau\to0$ and then $\delta\to0$
completes the proof.
\end{proof}

We now combine Lemma~\ref{L:2.3} with phase
control to prove Proposition~\ref{P:2.1}.
\begin{proof}[Proof of Proposition~\ref{P:2.1}]
Let us fix $g\in\cG$ and $\psi_0\in H^s$.
For sufficiently small $\tau>0$, let us set
$$
\psi_1=e^{ig/(2\tau)}\psi_0,\qquad
\psi_2=\cR_\tau(\psi_1,0),\qquad
\psi_3=e^{-i\varphi(1)/\tau}\psi_2.
$$
By Lemma~\ref{L:2.3}, $\psi_3\to\cL_g\psi_0$ in $H^s$
as $\tau\to0$. It remains to approximate $\psi_3$
by a controlled solution in time less than $2\tau$.
 
Let us fix such a $\tau$. By continuous dependence on the
initial state and continuity of multiplication by~$e^{-i\varphi(1)/\tau}$ in $H^s$, there is $\delta>0$
such that, whenever $\|v-\psi_1\|_s<\delta$, the
solution $\cR_t(v,0)$ exists on~$[0,\tau]$~and
$$
\big\|e^{-i\varphi(1)/\tau}\cR_\tau(v,0)
-\psi_3\big\|_s<\tau/2.
$$
By Theorem~\ref{T:1.6}, we can choose
$t_1\in(0,\tau/2)$ and $u_1\in\Theta(\psi_0,t_1)$
such that
$$
\|\cR_{t_1}(\psi_0,u_1)-\psi_1\|_s<\delta.
$$
After applying $u_1$, we let the solution evolve
without control for time $\tau$ and denote the
resulting state by~$h=\cR_\tau\bigl(\cR_{t_1}(\psi_0,u_1),0\bigr)$.
By the choice of $\delta$, we have
$$
\|e^{-i\varphi(1)/\tau}h-\psi_3\|_s<\tau/2.
$$
Applying Theorem~\ref{T:1.6} to $h$, we choose
$t_2\in(0,\tau/2)$ and $u_2\in\Theta(h,t_2)$ such that
$$
\|\cR_{t_2}(h,u_2)-e^{-i\varphi(1)/\tau}h\|_s<\tau/2.
$$
Concatenating $u_1$, zero control for time $\tau$,
and $u_2$, we obtain
$u_\tau\in\Theta(\psi_0,t_\tau)$, where
$t_\tau=t_1+\tau+t_2<2\tau$, such that
$$
\|\cR_{t_\tau}(\psi_0,u_\tau)-\psi_3\|_s<\tau.
$$
Since $\psi_3\to\cL_g\psi_0$ in $H^s$ as $\tau\to0$,
this proves the proposition.
\end{proof}

\subsection{Control of positive amplitudes}\label{SS:2.3}
 
In this subsection, we prove the following result
on the control of positive amplitudes and periodic phases.
 
\begin{proposition}\label{P:2.4}
Assume \textup{(\hyperlink{H1}{H$_1$})}.
Let $l\in\Z^d$ and
$a_0,a_1,\theta_0,\theta_1\in H^s(\T^d;\R)$ satisfy
$$
a_0,a_1>0,\qquad
\|a_0\|_{L^2}=\|a_1\|_{L^2}=\|\mathbf1\|_{L^2}.
$$
Then, for every $T>0$, equation~\eqref{E:0.1} is
approximately controllable from $a_0e^{i\theta_0}\phi_l$
to $a_1e^{i\theta_1}\phi_l$ at time $T$.
\end{proposition}
 
To prove this proposition, we use the following consequence
of Brenier's theorem on optimal transport~\cite{Bre},
in its periodic form~\cite{CE}, and of Caffarelli's
regularity theory~\cite{Caf}.
 
\begin{lemma}\label{L:2.5}
Let $s>d/2$ be an integer, and let
$\rho\in H^s(\T^d;\R)$ and
$\sigma\in C^\infty(\T^d;\R)$ be strictly positive
functions satisfying
$$
\int_{\T^d}\rho(x)\,\dd x
=\int_{\T^d}\sigma(x)\,\dd x=m.
$$
Then there is $g\in H^{s+2}(\T^d;\R)$ such that
$I+D^2g\ge cI$ for some $c>0$ and
\begin{equation}\label{E:2.11}
\sigma(x+\nabla g(x))\det(I+D^2g(x))=\rho(x),
\qquad x\in\T^d.
\end{equation}
\end{lemma}
 
\begin{proof}
Let us choose $0<\alpha<\min\{1,s-d/2\}$.
By Sobolev embedding, we have $\rho\in C^{0,\alpha}$. Applying the above-mentioned transport and regularity results in the form presented in
\cite[Theorem~4.3]{Loe}, to the
probability densities $\rho/m$ and $\sigma/m$, we find
 a convex function of the form
$$
\Phi(x)=\frac12|x|^2+g(x),
$$
with $g\in C^{2,\alpha}(\T^d;\R)$ satisfying~\eqref{E:2.11}.
Let us set $F=\Id+\nabla g$ and $A=I+D^2g$.
The convexity of $\Phi$ gives~$A=D^2\Phi\ge0$.
Since $\det A=\rho/\sigma(F)>0$, the matrix $A$
is positive definite, and by compactness,
$$
cI\le A(x)\le CI,\qquad x\in\T^d.
$$
 It remains to prove that $g\in H^{s+2}$.
Since $\sigma$ is smooth, the left-hand side
of~\eqref{E:2.11} is a smooth function
of~$x$, $\nabla g$, and $D^2g$.
Differentiating with respect to $g$ in a direction $h$
gives the linearized operator
$$
Lh=\sigma(F)(\det A)\tr(A^{-1}D^2h)
+\det(A)\lan\nabla\sigma(F),\nabla h\ran.
$$
Since $\sigma(F)\det A=\rho$, the coefficient matrix
of the second derivatives is $\rho A^{-1}$.
The bounds on $A$ and positivity of $\rho$ imply
$$
c_1I\le \rho(x)A(x)^{-1}\le C_1I,\qquad x\in\T^d.
$$
Thus~\eqref{E:2.11}
is uniformly elliptic at $g$.
We can therefore apply the nonlinear elliptic regularity
theorem \cite[Theorem~3.3.C]{TaylorPDE}, in its
second-order $L^2$ form with $r=\alpha$.
Since $g\in C^{2,\alpha}$ and $\rho\in H^s$,
it follows that~$g\in H^{s+2}(\T^d)$, which
completes the proof.
\end{proof}

\begin{proof}[Proof of Proposition~\ref{P:2.4}]
\emph{Step~1. Small-time controllability.}
Let us set $\psi_j=a_je^{i\theta_j}\phi_l$, $j=0,1$,
and fix~$\e>0$.
To apply Lemma~\ref{L:2.5} with a smooth target density,
let us choose a smooth positive function $b$ such that
$$
\|b\|_{L^2}=\|\mathbf1\|_{L^2},\qquad
\|(b-a_1)e^{i\theta_1}\phi_l\|_s<\e/2.
$$
Applying Lemma~\ref{L:2.5} with $\rho=a_0^2$ and
$\sigma=b^2$, we obtain $g\in\cG$ such that,
for $F=\Id+\nabla g$ and $P=F^{-1}$,
$$
b(F(x))^2\det DF(x)=a_0(x)^2,\qquad
(\det DP)^{1/2}a_0\circ P=b.
$$
It follows that
$$
\cL_g\psi_0=b e^{i\eta}\phi_l,\qquad
\eta=\theta_0\circ P+\lan l,P-\Id\ran.
$$
By Lemma~\ref{L:A.1}, $P-\Id$ is
periodic and belongs to $H^{s+1}$, and composition
by $P$ preserves $H^s$. Hence $\eta\in H^s(\T^d;\R)$
and
$$
e^{i(\theta_1-\eta)}\cL_g\psi_0
=b e^{i\theta_1}\phi_l.
$$
By Proposition~\ref{P:2.1}, we first approximate
$\cL_g\psi_0$ from $\psi_0$ in $H^s$.
We then apply Theorem~\ref{T:1.6} to the state
obtained, with phase $\theta_1-\eta$.
Since multiplication by $e^{i(\theta_1-\eta)}$
is bounded on $H^s$, concatenating these controls
gives a state arbitrarily close to
$be^{i\theta_1}\phi_l$.
Both controls can be chosen in arbitrarily small
time. Together with the choice of $b$, this proves
approximate controllability from $\psi_0$ to $\psi_1$
in arbitrarily small~time.
 
\emph{Step~2. Prescribed time.}
Let us fix $T>0$ and $\e>0$, and set $z=\phi_l$.
By Step~1, there are $t_1\in(0,T/3)$ and
$u_3\in\Theta(z,t_1)$ such that
$$
\|\cR_{t_1}(z,u_3)-\psi_1\|_s<\e/2.
$$
By Proposition~\ref{P:1.2}, we can choose $\delta_1>0$
so that, whenever $\|v-z\|_s<\delta_1$,
we have $u_3\in\Theta(v,t_1)$ and
$$
\|\cR_{t_1}(v,u_3)-\psi_1\|_s<\e.
$$
 
Let $\delta_0>0$ be given by Corollary~\ref{C:1.8}
with $T_0=T$ and accuracy $\delta_1$.
Using Step~1 again, we choose $t_0\in(0,T/3)$ and
$u_1\in\Theta(\psi_0,t_0)$ such that
$$
\|\cR_{t_0}(\psi_0,u_1)-z\|_s<\delta_0.
$$
Starting from this state, we apply Corollary~\ref{C:1.8}
to find a control $u_2$ on an interval of length
$T-t_0-t_1$ such that the final state is within
$\delta_1$-neighborhood of $z$.
Concatenating $u_1,u_2,u_3$, we obtain the required
approximate controllability property at time $T$.
\end{proof}

\section{Approximate controllability}\label{S:3}

In this section, we prove the following more general
version of the Main Theorem stated in the Introduction.

\begin{theorem}\label{T:3.1}
Let $s>d/2$ be an integer and assume
Condition~\textup{(\hyperlink{H1}{H$_1$})}.
For any nowhere vanishing initial and target states
$\psi_0,\psi_1\in H^s(\T^d;\C)\cap\bS$, any
$T>0$, and any $\e>0$, there are a control
$u\in L^2(0,T;\R^q)$ and a solution
$\psi\in C([0,T];H^s\cap\bS)$ of \eqref{E:0.1}
such that
$$
\psi(0)=\psi_0,\qquad \|\psi(T)-\psi_1\|_s<\e.
$$
\end{theorem}

When $\psi_0$ and $\psi_1$ have the same winding
vector, the result follows from
Proposition~\ref{P:2.4}.
Otherwise, we construct a short uncontrolled
trajectory whose initial and final states have
the winding vectors of $\psi_0$ and~$\psi_1$,
respectively.
We then use Proposition~\ref{P:2.4} to steer
the system from $\psi_0$ close to the initial state
of this trajectory, and from its final state
close to $\psi_1$.

\subsection{Proof of Theorem~\ref{T:3.1}}\label{SS:3.1}
 
The next lemma allows us to pass from one winding
vector to another along a short uncontrolled trajectory.
 
\begin{lemma}\label{L:3.2}
Let $l,m\in\Z^d$, $l\ne m$. There is $t_0>0$ such that,
for every $t\in(0,t_0]$, there are functions
$a_0,a_1,\theta_0,\theta_1\in C^\infty(\T^d;\R)$,
with $a_0,a_1>0$, such that the states
$$
\psi_a=a_0e^{i\theta_0}\phi_l,\qquad
\psi_b=a_1e^{i\theta_1}\phi_m
$$
belong to $\bS$ and satisfy $\cR_t(\psi_a,0)=\psi_b$.
\end{lemma}
 
This lemma is proved in Section~\ref{SS:3.2}.
 
\begin{proof}[Proof of Theorem~\ref{T:3.1}]
Let $\psi_0,\psi_1\in H^s\cap\bS$ be nowhere vanishing,
and fix $T>0$ and $\e>0$.
By Lemma~\ref{L:A.2}, there are $l,m\in\Z^d$
and $a_j,\theta_j\in H^s(\T^d;\R)$,
with $a_j>0$, such that
$$
\psi_0=a_0e^{i\theta_0}\phi_l,\qquad
\psi_1=a_1e^{i\theta_1}\phi_m.
$$
Since $\psi_0,\psi_1\in\bS$, we have
$\|a_0\|_{L^2}=\|a_1\|_{L^2}=\|\mathbf1\|_{L^2}$.
If $l=m$, the result follows directly from
Proposition~\ref{P:2.4}.
 
Let us assume that $l\ne m$.
By Lemma~\ref{L:3.2}, we can choose
$t\in(0,T)$ and states $\psi_a,\psi_b\in\bS$,
with the representations given there, such that
$\cR_t(\psi_a,0)=\psi_b$.
Let $r=(T-t)/2$. By Proposition~\ref{P:2.4},
we first choose $u_3\in\Theta(\psi_b,r)$ such that
$$
\|\cR_r(\psi_b,u_3)-\psi_1\|_s<\e/2.
$$
Let $v$ be the control obtained by using zero control
for time $t$, followed by $u_3$ for time $r$.
Then $v\in\Theta(\psi_a,t+r)$ and
$\cR_{t+r}(\psi_a,v)=\cR_r(\psi_b,u_3)$.
By Proposition~\ref{P:1.2}, there is $\eta>0$ such that,
whenever $\|h-\psi_a\|_s<\eta$, we have
$v\in\Theta(h,t+r)$ and
$$
\|\cR_{t+r}(h,v)-\psi_1\|_s<\e.
$$
Applying Proposition~\ref{P:2.4} again,
we choose $u_1\in\Theta(\psi_0,r)$ such that
$\|\cR_r(\psi_0,u_1)-\psi_a\|_s<\eta$.
Concatenating $u_1$ and $v$, we obtain
$u\in\Theta(\psi_0,T)$ satisfying
$$
\|\cR_T(\psi_0,u)-\psi_1\|_s<\e.
$$
This proves the theorem.
\end{proof}
 
\begin{remark}\label{R:3.3}
In the case $d=1$, the conclusion of
Theorem~\ref{T:3.1} remains true for any
target state $\psi_1\in H^s\cap\bS$, provided that the
initial state $\psi_0$ is nowhere vanishing. This follows
from the density of nowhere vanishing states
in~$H^s\cap\bS$. Indeed, let $\psi\in H^s\cap\bS$, and
let $f_n\in C^\infty(\T;\C)$ be such that $f_n\to\psi$
in $H^s$. Since $f_n(\T)$ is a smooth curve, it has empty
interior in $\C$, so we can find $c_n\in\C$ with
$|c_n|<1/n$ such that $f_n+c_n$ does not vanish.
The functions
$$
\psi_n=\frac{f_n+c_n}{\|f_n+c_n\|_{L^2}}
$$
are nowhere vanishing, belong to $H^s\cap\bS$, and
converge to~$\psi$ in $H^s$.

In the case $d\ge2$, the density property is no longer
true, since some zeros persist under small uniform
perturbations. Let us consider, for example, the state
$$
f(x)=(2\pi)^{-d/2}\big(\sin x_1+i\sin x_2\big),
$$
which belongs to $H^s\cap\bS$, and the disk
$$
D=\{x\in\T^d:x_1^2+x_2^2\le r^2,\ x_j=0
\text{ for }j\ge3\},\qquad 0<r<\pi.
$$
The function $f$ does not vanish on $\partial D$, and
$f/|f|$ has degree one there. If $g$ is a continuous
function sufficiently close to $f$ in $C^0$, then $g$ is
nonzero on $\partial D$, and $g/|g|$ also has degree one
on $\partial D$. Therefore, $g$ vanishes at some point
of~$D$. Since $H^s$ is continuously embedded in $C^0$,
we conclude that $f$ cannot be approximated in $H^s$ by
nowhere vanishing states. Let us note that this does not
imply that $f$ cannot be approximately reached by the
trajectories of~\eqref{E:0.1}; it only shows that the
density argument used for $d=1$ does not apply.

However, in any dimension, the conclusion of
Theorem~\ref{T:3.1} remains true for any
real-valued target state $\psi_1\in H^s\cap\bS$, even if
it has zeros, provided that $\psi_0$ is nowhere vanishing.
Indeed, for $\delta>0$, let us set
$$
\psi_{1,\delta}
=\frac{\psi_1+i\delta}{\big(1+(2\pi)^d\delta^2\big)^{1/2}}.
$$
Since $\psi_1$ is real-valued, we have
$|\psi_1+i\delta|^2=\psi_1^2+\delta^2\ge\delta^2$, so
$\psi_{1,\delta}$ is nowhere vanishing. Moreover, $\psi_{1,\delta}\in H^s\cap\bS$ and
$\psi_{1,\delta}\to\psi_1$ in $H^s$ as $\delta\to0$.
It remains to apply Theorem~\ref{T:3.1}
with the target $\psi_{1,\delta}$ for sufficiently
small $\delta>0$.
\end{remark}

\subsection{Construction of uncontrolled trajectories}\label{SS:3.2}
 
To prove Lemma~\ref{L:3.2}, we first construct
a free solution on the circle whose initial and
final states have winding numbers zero and one,
respectively; see Section~\ref{SS:3.2.1}.
We then replace the spatial variable in its initial
state by $\lan m-l,x\ran$ and multiply the resulting
function by $\phi_l$. This gives the initial state
for a free solution on $\T^d$.
In Section~\ref{SS:3.2.2}, we compare this
solution with the uncontrolled solution
of~\eqref{E:0.1} using Lemma~\ref{L:1.4}
and complete the proof.
 
\subsubsection{A one-dimensional construction}\label{SS:3.2.1}
 
\begin{lemma}\label{L:3.4}
There are constants $T_0>0$, $w_0\in(0,1]$, and $C>0$
such that, for every $w\in(0,w_0]$, there are functions
$h_w^0,h_w^1\in C^\infty(\T;\C)$ and
$\theta_w\in C^\infty(\T;\R)$ satisfying
\begin{align}
e^{iw^2T_0\partial_x^2}h_w^0&=h_w^1,
\label{E:3.1}\\
\|h_w^0-1\|_{L^\infty(\T)}&\le\frac34,
\label{E:3.2}\\
h_w^1(x)&=e^{ix}e^{i\theta_w(x)}.
\label{E:3.3}
\end{align}
Moreover,
\begin{equation}\label{E:3.4}
 \|h_w^0\|_{L^2(\T)}^2
=\|h_w^1\|_{L^2(\T)}^2=2\pi,
\qquad
\|h_w^0\|_{W(\T)}\le C.   
\end{equation}
\end{lemma}
 
For $d=1$, the winding vector in Lemma~\ref{L:A.2}
is an integer, which we call the \emph{winding number}.
It is the change of a continuous argument on
$[0,2\pi]$, divided by $2\pi$. 
By~\eqref{E:3.2}, $\Real h_w^0\ge1/4$,
so its argument can be chosen in~$(-\pi/2,\pi/2)$
and its winding number is zero.
For $h_w^1$, the argument $x+\theta_w(x)$ changes
by $2\pi$ over one period, since $\theta_w$ is
periodic. Thus its winding number is one.
 
The role of the parameter $w$ will be to make the transition
time $w^2T_0$ arbitrarily small while keeping
the Wiener norm of $h_w^0$ uniformly bounded.

Before proving Lemma~\ref{L:3.4}, we introduce
a periodization operator and establish some of its
elementary properties.
Let $\cS(\R)$ denote the Schwartz space of smooth
rapidly decreasing functions. For $w\in (0,1]$, we~define
$\cP_w:\cS(\R)\to C^\infty(\T;\C)$ by
$$
(\cP_wG)(x)=\sum_{k\in\Z}G\big((x+2\pi k)/w\big).
$$
For $G\in\cS(\R)$ and $h\in L^2(\T;\C)$, we denote
the Fourier transform of $G$ and the Fourier coefficients
of $h$ by
$$
\begin{aligned}
\hat G(\eta)&=\int_\R G(x)e^{-ix\eta}\,\dd x,
&&\eta\in\R,\\
\hat h(n)&=\frac1{2\pi}\int_\T h(x)e^{-inx}\,\dd x,
&&n\in\Z.
\end{aligned}
$$
This use of the same notation will cause no confusion.
 
\begin{lemma}\label{L:3.5}
For every $G\in\cS(\R)$ and $w\in (0,1]$, we have
\begin{align}
e^{it\partial_x^2}\cP_wG
&=\cP_w\big(e^{i(t/w^2)\partial_y^2}G\big),
\qquad t\in\R,\label{E:3.5}\\
\|\cP_wG\|_{W(\T)}
&\le\frac1{2\pi}\big(\|\hat G\|_{L^1(\R)}
+\|(\hat G)'\|_{L^1(\R)}\big),
\label{E:3.6}\\
\|\cP_wG\|_{L^\infty(\T)}
&\le\|G\|_{L^\infty(\R)}+C_Gw^2,
\label{E:3.7}
\end{align}
where $C_G>0$ depends only on $G$.
\end{lemma}
 
\begin{proof}
By the definition of $\cP_w$ and a change of variables,
$$
\widehat{\cP_wG}(n)
=\frac1{2\pi}\int_\R G(x/w)e^{-inx}\,\dd x
=\frac{w}{2\pi}\hat G(wn),\qquad n\in\Z.
$$
Therefore, both sides of \eqref{E:3.5} 
have Fourier coefficient
$(w/2\pi)e^{-itn^2}\hat G(wn)$ at $n$.
 
To prove \eqref{E:3.6}, we integrate
the inequality
$$
|\hat G(wn)|\le|\hat G(y)|+
\int_{wn}^{w(n+1)}|(\hat G)'(r)|\,\dd r,
\qquad wn\le y\le w(n+1)
$$
over this interval and sum over $n$. We obtain
$$
\|\cP_wG\|_{W(\T)}
=\frac1{2\pi}\sum_{n\in\Z}w|\hat G(wn)|
\le\frac1{2\pi}\big(
\|\hat G\|_{L^1}+w\|(\hat G)'\|_{L^1}\big),
$$
which proves~\eqref{E:3.6}
   since $w\le1$.
 
Finally, for $x\in[-\pi,\pi]$ and $k\ne0$, we have
$|x+2\pi k|\ge\pi|k|$. Since $G$ is rapidly decreasing,
$|G(y)|\le C|y|^{-2}$ for $y\ne0$, and hence
$$
\sum_{k\ne0}\left|G\big((x+2\pi k)/w\big)\right|
\le Cw^2\sum_{k\ne0}(\pi|k|)^{-2}.
$$
The term $k=0$ is bounded by $\|G\|_{L^\infty}$,
which proves \eqref{E:3.7}.
\end{proof}
 
\begin{proof}[Proof of Lemma~\ref{L:3.4}]
Let us choose a smooth nondecreasing function
$\chi:\R\to[0,2\pi]$ such that
$$
\chi(y)=0\quad\text{for }y\le-1,\qquad
\chi(y)=2\pi\quad\text{for }y\ge1,
$$
and set $G_1=e^{i\chi}-1$.
Then $G_1$ is smooth and supported in $[-1,1]$,
while $\mathbf1+G_1$ has modulus one and a phase
increasing from $0$ to $2\pi$. Next we set
$$
T_0=6,\qquad G_0=e^{-iT_0\partial_y^2}G_1.
$$
Thus $G_0\in\cS(\R)$ and
$e^{iT_0\partial_y^2}G_0=G_1$.
The kernel formula for $e^{it\partial_y^2}$ on $\R$
(e.g., see~\cite[(2.18)]{Tao})~gives
$$
G_0(x)=\frac{e^{i\pi/4}}{\sqrt{4\pi T_0}}
\int_\R e^{-i(x-y)^2/(4T_0)}G_1(y)\,\dd y.
$$
Since $G_1$ is supported in $[-1,1]$ and
$|G_1|\le2$, we have $\|G_1\|_{L^1(\R)}\le4$.
Therefore,
\begin{equation}\label{E:3.8}
\|G_0\|_{L^\infty(\R)}
\le\frac1{\sqrt{4\pi T_0}}\|G_1\|_{L^1(\R)}
\le\frac4{\sqrt{4\pi T_0}}<\frac12.
\end{equation}

We now rescale and periodize $G_0$ and $G_1$.
For $w\in (0,1]$, we set
$$
h_w^0=1+\cP_wG_0,\qquad
h_w^1=1+\cP_wG_1.
$$
Since $e^{it\partial_x^2}\mathbf1=\mathbf1$, the
fact that $e^{iT_0\partial_y^2}G_0=G_1$ and
\eqref{E:3.5} imply
\eqref{E:3.1}.
By~\eqref{E:3.7} and~\eqref{E:3.8}, we can choose~$w_0\in(0,1]$ sufficiently small that
\eqref{E:3.2} holds for $w\in(0,w_0]$.
The uniform bound on the Wiener norm in~\eqref{E:3.4} follows
from~\eqref{E:3.6}.
 
To verify~\eqref{E:3.3}, note that,
since $\operatorname{supp}G_1\subset[-1,1]$ and
$w\le1<\pi$, only the term $k=0$ contributes
to~$\cP_wG_1$ on $[-\pi,\pi]$. Hence
$h_w^1(x)=e^{i\chi(x/w)}$ on this interval.
Let us set
$$
\theta_w(x)=\chi(x/w)-x,\qquad -\pi\le x\le\pi.
$$
At $x=\pm\pi$, we have $\theta_w=\pi$,
$\theta_w'=-1$, and all higher derivatives vanish.
Thus $\theta_w$ extends to a function in~$C^\infty(\T;\R)$, and
\eqref{E:3.3} follows.
 
Finally, since $|h_w^1|=1$, unitarity of
$e^{it\partial_x^2}$ and~\eqref{E:3.1}
imply $\|h_w^0\|_{L^2(\T)}^2
=\|h_w^1\|_{L^2(\T)}^2=2\pi$.
\end{proof}

\subsubsection{Proof of Lemma~\ref{L:3.2}}\label{SS:3.2.2}

We first use Lemma~\ref{L:3.4} to construct
a short free trajectory on $\T^d$ with winding
vectors $l$ and $m$ at the initial and final times.
We then compare it with the uncontrolled solution
of~\eqref{E:0.1} starting from the same~state.

\emph{Step~1. Choice of the initial state.}
Let $k=m-l\ne0$. For $w\in(0,w_0]$, to be chosen
below, let~$h_w^0,h_w^1,\theta_w$ be given by
Lemma~\ref{L:3.4}.
Since $k=m-l$, equation~\eqref{E:3.3} gives
$$
\phi_l(x)h_w^1(\lan k,x\ran)
=\phi_l(x)e^{i\lan k,x\ran}
 e^{i\theta_w(\lan k,x\ran)}
=e^{i\theta_w(\lan k,x\ran)}\phi_m(x).
$$
We therefore choose
$$
\psi_a(x)=\phi_l(x)h_w^0(\lan k,x\ran).
$$
Since $\Real h_w^0>0$, the argument of
$h_w^0(\lan k,x\ran)$ with values in $(-\pi/2,\pi/2)$
defines a smooth real periodic function $\theta_0$.
Taking $a_0(x)=|h_w^0(\lan k,x\ran)|$, we obtain
$\psi_a=a_0e^{i\theta_0}\phi_l$, with
$a_0\in C^\infty(\T^d;\R)$ and $a_0>0$.
 
To verify that $\psi_a\in\bS$, note that, for every
integrable function $H:\T\to\C$,
$$
\int_{\T^d}H(\lan k,x\ran)\,\dd x
=(2\pi)^{d-1}\int_\T H(y)\,\dd y.
$$
Applying this identity to $H=|h_w^0|^2$ and using
\eqref{E:3.4}, we get
$$
\|\psi_a\|_{L^2}^2
=\frac1{2\pi}\|h_w^0\|_{L^2(\T)}^2=1.
$$
 
\emph{Step~2. Computation of $e^{it\Delta}\psi_a$.}
For $h\in C^\infty(\T;\C)$, Fourier expansion gives
\begin{equation}\label{E:3.9}
\phi_l(x)h(\lan k,x\ran)
=\sum_{n\in\Z}\hat h(n)\phi_{l+nk}(x).
\end{equation}
The operator $e^{it\Delta}$ multiplies each term
by $e^{-it|l+nk|^2}$. Since
$$    
|l+nk|^2=|l|^2+2n\lan l,k\ran+n^2|k|^2,
$$
we obtain
\begin{equation}\label{E:3.10}
e^{it\Delta}\big(\phi_l(x)h(\lan k,x\ran)\big)
=e^{-it|l|^2}\phi_l(x)
\big(e^{it|k|^2\partial_y^2}h\big)
(\lan k,x\ran-2t\lan l,k\ran).
\end{equation}
To make the time in the one-dimensional evolution
equal to $w^2T_0$, let us set
$$
t_w=\frac{w^2T_0}{|k|^2},\qquad
\delta_w=2t_w\lan l,k\ran.
$$
By \eqref{E:3.1},  
\eqref{E:3.3}, and \eqref{E:3.10}, the state
$\psi_f=e^{it_w\Delta}\psi_a$ satisfies
$$
\psi_f(x)
=e^{-it_w|l|^2}\phi_l(x)h_w^1(\lan k,x\ran-\delta_w)
=e^{i\eta_w(x)}\phi_m(x),
$$
where
$$
\eta_w(x)=\theta_w(\lan k,x\ran-\delta_w)
-t_w|l|^2-\delta_w.
$$
Since $\theta_w$ is periodic and $k\in\Z^d$,
we have $\eta_w\in C^\infty(\T^d;\R)$.
In particular, $|\psi_f|=(2\pi)^{-d/2}$.
 
\emph{Step~3. Comparison with the uncontrolled equation.} The Sobolev norms of $\psi_a$ may grow as $w\to0$.
To show that the final state of the uncontrolled
solution has winding vector $m$, it is enough
to prove that it is sufficiently close to $\psi_f$
in the uniform norm.
To apply Lemma~\ref{L:1.4}, we need
a bound on the Wiener norm of $\psi_a$
independent of $w$.
Since $k\ne0$, the map $n\mapsto l+nk$ from $\Z$
to $\Z^d$ is injective.
Thus the Fourier expansion~\eqref{E:3.9}
gives
$$
\|\psi_a\|_{W(\T^d)}
=(2\pi)^{-d/2}\|h_w^0\|_{W(\T)}\le N,
$$
where $N=(2\pi)^{-d/2}C$ and $C$ is the constant
in \eqref{E:3.4}. Thus $N$ is independent
of $w\in(0,w_0]$.
Let $T_N,C_N$ be given by
Lemma~\ref{L:1.4}, and choose $t_0>0$ so small that
$$
|k|\sqrt{t_0/T_0}\le w_0,\qquad
t_0<T_N,\qquad
(2\pi)^{d/2}C_Nt_0<\frac14.
$$
For any $t\in(0,t_0]$, let us take
$w=|k|\sqrt{t/T_0}$. Then $w\in(0,w_0]$ and $t_w=t$.
By Lemmas~\ref{L:1.3} and~\ref{L:1.4}, $\psi_b=\cR_t(\psi_a,0)$ is well defined, smooth,
and satisfies
$$
\|\psi_b-\psi_f\|_{W(\T^d)}\le C_Nt.
$$
Since $|\psi_f|=(2\pi)^{-d/2}$,
this estimate also controls the quotient
$\psi_b/\psi_f$. Indeed, using
$\|\cdot\|_{L^\infty}\le\|\cdot\|_W$ (see~\eqref{E:1.1}), we obtain
$$
\left\|\frac{\psi_b}{\psi_f}-1\right\|_{L^\infty}
=(2\pi)^{d/2}\|\psi_b-\psi_f\|_{L^\infty}
\le(2\pi)^{d/2}C_Nt<\frac14.
$$
Thus $\Real(\psi_b/\psi_f)>3/4$, so its argument
can be chosen smoothly and periodically with values
in $(-\pi/2,\pi/2)$. We can therefore write
$$
\frac{\psi_b}{\psi_f}=a_1e^{i\beta_w},
\qquad
a_1=\left|\frac{\psi_b}{\psi_f}\right|>0,
$$
where $a_1,\beta_w\in C^\infty(\T^d;\R)$.
It follows that
$$
\psi_b=a_1e^{i(\eta_w+\beta_w)}\phi_m.
$$
Taking $\theta_1=\eta_w+\beta_w$, we complete the proof.

\appendix
\refstepcounter{section}
\section*{Appendix}
\addtocontents{toc}{\protect\setcounter{tocdepth}{1}}
\addcontentsline{toc}{section}{Appendix}
\label{S:A}

\subsection{Proof of Corollary~\ref{C:1.8}}\label{SS:A.1}
 
We consider equation~\eqref{E:0.1} with smooth real-valued potentials
$V,Q_1,\dots,Q_q$ and assume  Condition~\textup{(\hyperlink{H1}{H$_1$})}. Let~$z=\phi_l$, fix any 
 $T_0>0$ and $\e>0$, and denote
\begin{equation}\label{E:A.1}
V_z=V+|l|^2+\kappa|z|^{2p}.
\end{equation}
We adapt the fixed-time argument of \cite[Theorem~3.3, Step~4]{Ner21}, alternating short intervals with zero control and phase controls. Multiplication by $e^{ihV_z}$ cancels the
first-order change of the uncontrolled solution
starting from $z$ over time $h$.
The remaining error is of order $h^2$, which allows
us to repeat the construction while keeping the
accumulated error small.
 The proof is divided into five steps.
 
\emph{Step~1. Short-time estimates.} Let us choose $h_0>0$ so small that $U(t):=\cR_t(z,0)$ exists on $[0,h_0]$ and
$$
\sup_{0\le t\le h_0}\|U(t)-z\|_s<\frac\e2.
$$
By Proposition~\ref{P:1.2}, there is $r_0>0$ such that, whenever $\|w-z\|_s<r_0$, the solution $\cR_t(w,0)$ exists on~$[0,h_0]$~and  
$$
\sup_{0\le t\le h_0}\|\cR_t(w,0)-U(t)\|_s<\frac\e2.
$$
Consequently,
$$
\sup_{0\le t\le h_0}\|\cR_t(w,0)-z\|_s<\e.
$$
We shall use this estimate in Step~5. It also gives
$$
\sup_{0\le t\le h_0}\|\cR_t(w,0)\|_s\le\|z\|_s+\e,
 \quad \|w-z\|_s<r_0.
$$
Thus the solutions $\cR_t(w,0)$ are uniformly bounded in $H^s$ on $[0,h_0]$.
 
\emph{Step~2. Estimate after phase multiplication.}
We claim that there are constants $C,K\ge0$ such that
\begin{equation}\label{E:A.2}
\big\|e^{ihV_z}\cR_h(w,0)-z\big\|_s
\le e^{Ch}\|w-z\|_s+Kh^2,
\qquad 0\le h\le h_0
\end{equation}
for every $w\in H^s$ with $\|w-z\|_s<r_0$. To prove this, let us set
$$
F(f)=Vf+\kappa|f|^{2p}f,\qquad f\in H^s.
$$
Since $z$ is smooth, Lemma~\ref{L:1.3} gives $U\in C([0,h_0];H^{s+4})$. The equation therefore implies
$$
U'=i\Delta U-iF(U)\in C([0,h_0];H^{s+2}).
$$
Differentiating one more time, we obtain
$$
U''=i\Delta U'-iF'(U)U'\in C([0,h_0];H^s),
$$
where $F'$ denotes the real Fr\'echet derivative. Thus $y(t):=e^{itV_z}U(t)$ belongs to $C^2([0,h_0];H^s)$. Using $-\Delta z=|l|^2z$ and~\eqref{E:A.1}, we get
$$
U'(0)=-iV_z z,\qquad y(0)=z,\qquad
y'(0)=iV_z z+U'(0)=0.
$$
This implies
$$
y(h)-z
=\int_0^h(h-r)y''(r)\,\dd r,
$$
and consequently,
\begin{equation}\label{E:A.3}
    \|e^{ihV_z}\cR_h(z,0)-z\|_s\le Kh^2,
\qquad 0\le h\le h_0,
\end{equation}
where
$$
K=\frac12\max_{0\le t\le h_0}\|y''(t)\|_s<\infty.
$$

We next compare the solutions starting from $w$ and $z$. Since $s>d/2$, there is $C_0\ge0$ such that
\begin{equation}\label{E:A.4}
\|F(f)-F(g)\|_s\le C_0\|f-g\|_s
\end{equation}
for all $f,g\in H^s$ with $\|f\|_s,\|g\|_s\le\|z\|_s+\e$. By Step~1, this estimate applies to $f=\cR_t(w,0)$ and $g=U(t)$ for $0\le t\le h_0$. Subtracting the Duhamel formulas for the solutions 
$\cR_t(w,0)$ and $U(t)$ gives
$$
\cR_t(w,0)-U(t)
=e^{it\Delta}(w-z)
-i\int_0^t e^{i(t-r)\Delta}
\bigl(F(\cR_r(w,0))-F(U(r))\bigr)\,\dd r.
$$
Using the isometry of $e^{it\Delta}$ and \eqref{E:A.4}, we obtain
$$
\begin{aligned}
\|\cR_t(w,0)-U(t)\|_s
&\le \|w-z\|_s
+\int_0^t\|F(\cR_r(w,0))-F(U(r))\|_s\,\dd r\\
&\le \|w-z\|_s
+C_0\int_0^t\|\cR_r(w,0)-U(r)\|_s\,\dd r.
\end{aligned}
$$
By the Gronwall inequality,
\begin{equation}\label{E:A.5}
\|\cR_t(w,0)-U(t)\|_s
\le e^{C_0t}\|w-z\|_s,
\qquad 0\le t\le h_0.
\end{equation}
 
We also need an $H^s$ estimate for the phase multiplication. Since $V_z$ is smooth, there is $C_1\ge0$ such that $\|V_zf\|_s\le C_1\|f\|_s$ for all $f\in H^s$. The identity
$$
e^{itV_z}f=f+i\int_0^tV_ze^{irV_z}f\,\dd r
$$
and Gronwall's inequality give
$$
\|e^{itV_z}f\|_s\le e^{C_1t}\|f\|_s,\qquad 0\le t\le h_0.
$$
Combining this with \eqref{E:A.3} and \eqref{E:A.5}, we obtain
$$
\begin{aligned}
\|e^{ihV_z}\cR_h(w,0)-z\|_s
&\le e^{C_1h}\|\cR_h(w,0)-\cR_h(z,0)\|_s+Kh^2\\
&\le e^{(C_0+C_1)h}\|w-z\|_s+Kh^2.
\end{aligned}
$$
This proves~\eqref{E:A.2} with $C=C_0+C_1$. The constants $h_0,r_0,C,K$ do not depend on $T$.
 
\emph{Step~3. Construction of one control.} Let us fix $h\in(0,h_0]$, $\eta>0$, and $w\in H^s$ with
$\|w-z\|_s<r_0$. Starting from $w$, we use zero control for time $h$ and set $\tilde w=\cR_h(w,0)$.
 By Theorem~\ref{T:1.6}, applied to $\tilde w$ with phase $hV_z$, there are a time $0<\tau<\eta$ and a control $u\in\Theta(\tilde w,\tau)$ such that
$$
\big\|\cR_\tau(\tilde w,u)-e^{ihV_z}\tilde w\big\|_s<h^2.
$$
Combining this  with~\eqref{E:A.2}, we obtain
$$
\begin{aligned}
\|\cR_\tau(\tilde w,u)-z\|_s
&\le \|\cR_\tau(\tilde w,u)-e^{ihV_z}\tilde w\|_s
   +\|e^{ihV_z}\tilde w-z\|_s\\
&\le e^{Ch}\|w-z\|_s+(K+1)h^2.
\end{aligned}
$$
We shall choose $\eta$ in the next step.

\emph{Step~4. Iteration and choice of the parameters.}
For $h\in(0,h_0]$ and $v\in H^s$ with $\|v-z\|_s<r_0$,
we define a sequence $\{v_k\}$ as follows. We set $v_0=v$, and
given $v_k$ with $\|v_k-z\|_s<r_0$, we define
$$
v_{k+1}=\cR_{\tau_k}\bigl(\cR_h(v_k,0),u_k\bigr),
$$
where the time $\tau_k$ and control $u_k$ are given by
Step~3 with $w=v_k$. Then
$$
\|v_{k+1}-z\|_s
\le e^{Ch}\|v_k-z\|_s+(K+1)h^2.
$$
Iterating this estimate, we obtain
\begin{equation}\label{E:A.6}
\|v_k-z\|_s
\le e^{Ckh}\|v-z\|_s
  +(K+1)h^2\sum_{j=0}^{k-1}e^{Cjh}
\le e^{CT_0}\bigl(\|v-z\|_s+(K+1)T_0h\bigr),
\end{equation}
provided $kh\le T_0$ and $\|v_j-z\|_s<r_0$ for
$j=0,\ldots,k-1$.
 
Let us choose $\delta>0$ such that $e^{CT_0}\delta<r_0$. Next we fix $T\in(0,T_0]$
and $v\in H^s$ with $\|v-z\|_s<\delta$, and choose an
integer $N\ge2$ so large that $h=T/N$ satisfies
$$
h\le h_0,\qquad
e^{CT_0}\bigl(\delta+(K+1)T_0h\bigr)<r_0.
$$

We use the construction above with $\eta=h/(2N)$ in Step~3.
Estimate~\eqref{E:A.6}
and the choices of $\delta$ and $h$ imply by induction that
all $N-1$ steps can be carried out so that
$$
\|v_k-z\|_s<r_0,\qquad k=0,\ldots,N-1,
$$
and the corresponding times satisfy
$0<\tau_k<h/(2N)$ for $k=0,\ldots,N-2$. The state $v_{N-1}$ is reached at~time
$$
t_{N-1}:=(N-1)h+\sum_{k=0}^{N-2}\tau_k.
$$
In the final step, we will check that $t_{N-1}<T$ and extend the
control by zero from time $t_{N-1}$ up to~$T$.

\par\bigskip
\newsavebox{\returnfigbox}%
\noindent\begin{minipage}{\linewidth}
\centering
\captionsetup{type=figure,skip=7pt}
\usetikzlibrary{decorations.markings}
\sbox{\returnfigbox}{%
\begin{tikzpicture}[x=1.08cm,y=1.35cm,
  every node/.style={font=\small,inner sep=1.5pt},
  drift/.style={draw=black!83,line width=.92pt},
  kick/.style={draw=darkblue,line width=1.22pt},
  ideal/.style={draw=teal!64!black,line width=.95pt,
                dash pattern=on 3.4pt off 2.1pt},
  separation/.style={draw=black!58,line width=.58pt,densely dotted},
  direction/.style={postaction={decorate},
    decoration={markings,mark=at position #1 with
      {\arrow{Stealth[length=4.8pt,width=3.9pt]}}}},
  state/.style={circle,fill=darkblue,draw=white,
                line width=.45pt,inner sep=2pt},
  idealpoint/.style={circle,fill=white,draw=teal!64!black,
                     line width=.85pt,inner sep=1.7pt}]
\coordinate (center) at (2.02,1.12);
\coordinate (start) at (1.15,2.30);
\coordinate (driftend) at (11.65,2.35);
\coordinate (idealend) at (3.50,1.55);
\coordinate (actualend) at (3.50,.66);
 
\filldraw[fill=black!3,draw=black!45,line width=.65pt,
          dash pattern=on 3pt off 2.5pt]
  (center) ellipse [x radius=1.90,y radius=1.50];
\node[text=black!65] at (1.65,.25) {$B_{H^s}(z,r_0)$};
 
\draw[separation] (idealend) -- (actualend);
\draw[ideal,direction=.52] (driftend)
  .. controls (9.45,2.55) and (6.45,1.55) .. (idealend);
\draw[drift,direction=.53] (start)
  .. controls (3.65,4.10) and (8.75,3.80) .. (driftend);
\draw[kick,direction=.54] (driftend)
  .. controls (11.30,-.07) and (6.80,-.48) .. (actualend);
 
\fill[black!80] (center) circle (1.45pt);
\node[below left=2pt] at (center) {$z$};
\fill[black!85] (start) circle (2pt);
\node[above left=2pt] at (start) {$v_k$};
\fill[black!82] (driftend) circle (1.45pt);
\node[right=4pt] at (driftend) {$\tilde v_k$};
\node[idealpoint] at (idealend) {};
\node[above left=2pt,text=teal!64!black] at (idealend)
  {$e^{ihV_z}\tilde v_k$};
\node[state] at (actualend) {};
\node[below left=2pt] at (actualend) {$v_{k+1}$};
 
\node[above=3pt] at (6.30,3.60) {time $h$};
\node[above=3pt,text=teal!64!black] at (7.75,2.00) {$e^{ihV_z}$};
\node[below=3pt,text=darkblue] at (8.25,.16) {$u_k$, time $\tau_k$};
\draw[black!45,line width=.4pt] (3.63,1.10) -- (4.20,1.10);
\node[anchor=west] at (4.28,1.10)
  {$\|v_{k+1}-e^{ihV_z}\tilde v_k\|_s<h^2$};
 
\begin{scope}[yshift=-.70cm]
  \draw[drift] (.75,0) -- (1.38,0);
  \node[right] at (1.47,0) {zero control};
  \draw[kick,direction=.55] (4.35,0) -- (4.98,0);
  \node[right] at (5.07,0) {phase control};
  \draw[ideal] (8.10,0) -- (8.73,0);
  \node[right] at (8.82,0) {phase multiplication};
\end{scope}
\end{tikzpicture}}%
\vbox to 1.06\dimexpr\ht\returnfigbox+\dp\returnfigbox\relax{%
  \vss\hbox to\linewidth{\hss\scalebox{0.96}{\usebox{\returnfigbox}}\hss}\vss}%
\bigskip
\caption{One step of the return construction. From $v_k$, we follow
 the uncontrolled flow to $\tilde v_k=\cR_h(v_k,0)$ and apply
 the phase control $u_k$ to reach $v_{k+1}$. The dotted segment
 shows the error relative to $e^{ihV_z}\tilde v_k$.}
\label{F:return-step}
\end{minipage}
\par\bigskip\bigskip

\emph{Step~5. Completion of the proof.}
Since $T=Nh$ and $0<\tau_k<h/(2N)$, the remaining time
satisfies
$$
r:=T-t_{N-1}=h-\sum_{k=0}^{N-2}\tau_k,
\qquad \frac h2<r<h\le h_0.
$$
In particular, $t_{N-1}<T$. Extending the control by zero on $[t_{N-1},T]$, we obtain
$u\in L^2(0,T;\R^q)$.
 
Since $\|v_{N-1}-z\|_s<r_0$ and $r\le h_0$, Step~1
shows that the solution with initial data $v_{N-1}$ and
zero control exists on $[0,r]$. Thus the constructed solution extends to time $T$, so
$u\in\Theta(v,T)$, and the estimate in Step~1 gives
$$
\|\cR_T(v,u)-z\|_s
=\|\cR_r(v_{N-1},0)-z\|_s<\e,
$$
which completes the proof.

\subsection{Inverse maps and composition}\label{SS:A.2}
 
The following lemma gives the properties of the
transport maps and operators used in
Section~\ref{S:2}.
 
\begin{lemma}\label{L:A.1}
Let $s>d/2$ be an integer and let $g\in\cG$.
Then the following properties hold.
\begin{itemize}
\item[(i)] For every $\sigma\in[0,1]$, the maps
$F_\sigma$ and $f_\sigma$ are $C^1$ diffeomorphisms
of $\R^d$ and $\T^d$, respectively.
 
\item[(ii)] Let us set $P_\sigma=F_\sigma^{-1}$.
The maps $P_\sigma$ and $f_\sigma^{-1}$ are jointly
$C^1$ in $\sigma$ and the spatial variable.
The function $P_\sigma-\Id$ is $2\pi\Z^d$-periodic
and depends continuously on $\sigma$ in
$H^{s+1}(\T^d;\R^d)$.
For $k=s,s+1$ and every fixed $h\in H^k$,
the map $\sigma\mapsto h\circ P_\sigma$ is continuous
in $H^k$, and
$$
\sup_{\sigma\in[0,1]}\|h\circ P_\sigma\|_k
\le C\|h\|_k.
$$
 
\item[(iii)] For every fixed $h\in H^s$, the map
$\sigma\mapsto\cL_{\sigma g}h$ is continuous
in $H^s$, and
$$
\sup_{\sigma\in[0,1]}\|\cL_{\sigma g}h\|_s
\le C\|h\|_s.
$$
\end{itemize}
Here $C>0$ depends only on $g$ and $s$.
\end{lemma}
 
\begin{proof}
\emph{(i)} By the definition of $\cG$ and compactness,
there is $c_0>0$ such that $I+D^2g(x)\ge c_0I$ for
every~$x\in\R^d$. Thus, for $\sigma\in[0,1]$,
$$
DF_\sigma(x)
=(1-\sigma)I+\sigma(I+D^2g(x))
\ge(1-\sigma+\sigma c_0)I\ge cI,
$$
where $c=\min\{1,c_0\}$. It follows that
$$
\lan F_\sigma(x)-F_\sigma(y),x-y\ran
\ge c|x-y|^2,
$$
and the Cauchy--Schwarz inequality gives
$$
|F_\sigma(x)-F_\sigma(y)|\ge c|x-y|.
$$
Hence $F_\sigma$ is injective and has closed range.
Its range is also open by the inverse mapping
theorem. Since~$\R^d$ is connected, $F_\sigma$
is surjective and hence a $C^1$ diffeomorphism. By~\eqref{E:2.1}, $F_\sigma$ induces
a diffeomorphism $f_\sigma$ of~$\T^d$.

\emph{(ii)} Let $P_\sigma=F_\sigma^{-1}$.
By~\eqref{E:2.1}, we have
$$
P_\sigma(x+2\pi n)=P_\sigma(x)+2\pi n,
\qquad n\in\Z^d.
$$
Thus $r_\sigma=P_\sigma-\Id$ is periodic and
$f_\sigma^{-1}([x])=[x+r_\sigma(x)]$.
Applying the inverse mapping theorem to~$(\sigma,x)\mapsto(\sigma,F_\sigma(x))$ gives joint
$C^1$ dependence of the inverses on $\sigma$ and
the spatial variable. If $g$ is smooth, the maps
and their inverses are smooth.
 
Since $F_\sigma-\Id=\sigma\nabla g$ is continuous
in $H^{s+1}$ and $s+1>d/2+1$, the inversion and
composition results of
\cite[Lemmas~2.7--2.8 and Remark~1.3]{IKT}
apply to the induced torus maps. Together with the equality~$P_\sigma-\Id=-\sigma\nabla g\circ P_\sigma$,
they give the continuity of $P_\sigma-\Id$ in
$H^{s+1}$.
 
The same composition result applies to
functions in $H^k$, $k=s,s+1$, and gives continuity
of $\sigma\mapsto h\circ P_\sigma$ for every fixed
$h\in H^k$. The corresponding bound is uniform in $\sigma\in[0,1]$,
since $P_\sigma-\Id$ is continuous in~$H^{s+1}$ and
$\det DP_\sigma$ has a positive minimum on
$[0,1]\times\T^d$.
 
\emph{(iii)} Since $DP_\sigma$ is continuous in $H^s$
and its determinant is bounded away from zero,
the algebra property of $H^s$  implies that
$(\det DP_\sigma)^{1/2}$ is continuous in $H^s$.
The identity
$$
\cL_{\sigma g}h=(\det DP_\sigma)^{1/2}(h\circ P_\sigma)
$$
and multiplication in $H^s$ give the continuity
in~\textup{(iii)}. Its uniform bound follows from
\textup{(ii)} and boundedness of~$(\det DP_\sigma)^{1/2}$ in $H^s$ for $\sigma\in[0,1]$.
\end{proof}

\subsection{Representation of nowhere vanishing states}\label{SS:A.3}
 
\begin{lemma}\label{L:A.2}
Let $s>d/2$ be an integer. For every nowhere vanishing
$\psi\in H^s(\T^d;\C)$, there is a unique vector
$l\in\Z^d$, called the winding vector of $\psi$,
such that $\psi=ae^{i\theta}\phi_l$ for some
$a,\theta\in H^s(\T^d;\R)$ with $a>0$.
\end{lemma}
\begin{proof}
Since $\psi$ is continuous and nowhere vanishing,
$\min_{x\in\T^d}|\psi(x)|>0$. Let us consider the functions
$$
a=(2\pi)^{d/2}|\psi|,\qquad b=\frac{\psi}{|\psi|}.
$$
Then $a\in H^s(\T^d;\R)$ with $a>0$, and
$b\in H^s(\T^d;\C)$ with $|b|=1$. 
 
We extend $b$ periodically to $\R^d$ and choose~$\alpha\in\R$ such that $e^{i\alpha}=b(0)$.
For each $n\ge1$, let 
$$
B_n=\{x\in\R^d:|x|<n\}.
$$
Applying the  Sobolev lifting property in~\cite[Corollary~1]{BBM00} to~$b|_{B_n}$,
we obtain $\vartheta_n\in H^s(B_n;\R)$ such that
$b=e^{i\vartheta_n}$ on $B_n$.
By Sobolev embedding, $\vartheta_n$ is continuous.
Adding an integer multiple of $2\pi$, we may assume
that $\vartheta_n(0)=\alpha$.
For~$m>n$, the difference $\vartheta_m-\vartheta_n$
is continuous on $B_n$, takes values in~$2\pi\Z$,
and vanishes at the origin. Since $B_n$ is connected,
$\vartheta_m=\vartheta_n$ on $B_n$.
Thus these functions define~$\vartheta\in H^s_{\mathrm{loc}}(\R^d;\R)$ satisfying
$b=e^{i\vartheta}$.

For each $j=1,\dots,d$, periodicity of $b$ implies
$$
\vartheta(x+2\pi e_j)-\vartheta(x)\in2\pi\Z,
\qquad x\in\R^d.
$$
By continuity, this difference is constant, equal to
$2\pi l_j$ for some $l_j\in\Z$.
Let $l=(l_1,\dots,l_d)$ and set
$$
\theta(x)=\vartheta(x)-\lan l,x\ran.
$$
Then $\theta$ is $2\pi\Z^d$-periodic and belongs to
$H^s_{\mathrm{loc}}(\R^d;\R)$, so
$\theta\in H^s(\T^d;\R)$.
Finally, $\psi=|\psi|e^{i\vartheta}=ae^{i\theta}\phi_l$.
The~uniqueness of $l$ follows from the fact that
any two continuous arguments of $b$ on $\R^d$
differ by a constant~in~$2\pi\Z$.
\end{proof}

\section*{Declaration of generative AI use}
 
ChatGPT suggested the construction of the special trajectory in Lemma \ref{L:3.4}. The remaining arguments are the author’s.  ChatGPT has been used to  check the proofs, to proofread and polish the text, and to assist in preparing the figure in the proof of Corollary 1.8. The author takes full responsibility for the~paper.

\section*{Acknowledgments}
The author thanks Eugenio Pozzoli for   discussions
of the  methods of~\cite{BPexamples,BP,BCP}
during his visit to NYU Shanghai in June 2026.
This work was partially supported by the National
Natural Science Foundation of China
(Grant No.~12571156) and the NYU Shanghai Boost Fund.

\end{document}